# Method for Generating Distributions and Classes of Probability Distributions: the Univariate Case

Cícero Carlos Ramos de Brito

*Statistics and Informatics Department, Universidade Federal Rural de Pernambuco, Recife, Brazil.*

Address: Av. Professor Cláudio Selva, 178, Dois Irmãos, Recife, Pernambuco, Brazil. CEP: 52171-260

cicerocarlosbrito@yahoo.com.br

Leandro Chaves Rêgo

*Department of Statistics, Universidade Federal de Pernambuco, Recife, Brazil*

leandro@de.ufpe.br

Wilson Rosa de Oliveira

*Statistics and Informatics Department, Universidade Federal Rural de Pernambuco, Recife, Brazil.*

wilson.rosa@gmail.com

# Method to Generate Distributions and Classes of Probability Distributions: the Univariate Case

In this work, we present a method to generate probability distributions and classes of probability distributions, which broadens a process of probability distribution construction. In this method, distribution classes are built from pre-defined monotonic functions and from known distributions. With the use of this method, we can obtain different classes of probability distributions described in literature. Beside these results, we could obtain results on the support and nature of the generated distributions.

KEYWORDS: Probabilistic distributions generating method, Distribution classes generating functions, probabilistic distributions classes, probabilistic distribution support.

## 1. Introduction

The amount of data available for analysis is growing increasingly faster, requiring new probabilistic distributions to better describe each phenomenon or experiment studied. Computer based tools allow the use of more complex distributions with a larger number of parameters to better study sizeable masses of data.

The literature in the field describes several generalizations and extensions of symmetric, asymmetric, discrete and continuous distributions. The relevance of these new models is that, according to the situation, each one of them can better fit the mass of data. Table 1 presents several classes of distributions described in literature, their nomenclature and the title of the work where they have been presented.

TABLE 1: Some classes of distributions described in literature

| Distribution classes | Nomenclature |
|---|---|
| $F(x) = G^a(x)$, where $a > 0$ | Exponentiated generalized defined by Mudholkar *et al.* (1995) |
| $F(x) = \frac{1}{B(a,b)}\int_0^{G(x)} t^{a-1}(1-t)^{b-1}dt$, where $a > 0$, $b > 0$ and $0 < t < 1$. | beta1 generalized defined by Eugene *et al.* (2002) |
| $F(x) = \frac{1}{B(a,b)} \int_0^{G(x)} t^{a-1}(1+t)^{-(a+b)} dt$, where $a > 0$, $b > 0$ and $t > 0$ | Beta3 generalized defined by Thair and Nadarajah (2013) |
| $F(x) = \frac{1}{B(a,b)} \int_0^{G^c(x)} t^{a-1}(1-t)^{b-1} dt$, where $a > 0$, $b > 0$, c> 0 and $0 < t < 1$ | Mc1 generalized defined by McDonald (1984) |
| $F(x) = \frac{1}{B(a,b)} \int_0^{G^c(x)} t^{a-1}(1+t)^{-(a+b)} dt$, where $a > 0$, $b > 0$, c> 0 and $t > 0$ | Mc3 generalized defined by Thair and Nadarajah (2013) |
| $F(x) = 1 - \left(1 - G^a(x)\right)^b$ | Kumaraswamy $G_1$ defined by Cordeiro and Castro (2011) |
| $F(x) = 1 - \left(1 - \left(1 - G(x)\right)^a\right)^b$, where $a > 0$ and $b > 0$ | Kumaraswamy type 2 defined by Thair and Nadarajah (2013) |
| $F(x) = \frac{G(x)}{G(x) + b\left(1 - G(x)\right)}$, where $b > 0$ | Marshall-Olkin defined by Marshall and Olkin (1997) |
| $F(x) = 1 - \left(\frac{b\left(1-G(x)\right)}{G(x) + b\left(1 - G(x)\right)}\right)^\theta$, where $b > 0$ and $\theta > 0$ | Marshall-Olkin $G_1$ defined by Jayakumar and Mathew (2008) |
| $F(x) = \left(\frac{G(x)}{G(x) + b\left(1 - G(x)\right)}\right)^\theta$, where $b > 0$ and $\theta > 0$ | Marshall-Olkin $G_1$ defined by Thair and Nadarajah (2013) |
| $$F(x) = \frac{\beta^\alpha}{\Gamma(\alpha)} \int_0^{-ln\left(1-G(x)\right)} t^{\alpha-1} e^{-\beta t} dt$$ | Gamma-Generated defined by Zografos and Balakrishnan (2009) |
| $$F(x) = 1 - \frac{\beta^\alpha}{\Gamma(\alpha)} \int_0^{-ln\left(G(x)\right)} t^{\alpha-1} e^{-\beta t} dt$$ | Gamma-Generated defined by Cordeiro (2013) |
| $F(x) = 1 - \frac{C\left(\theta e^{-\alpha H(x)}\right)}{C(\theta)}$, where $x > 0$, $\theta > 0$ and $C(\theta) = \sum_{n=1}^{\infty} a_n \theta^n$ | Extended Weibull distribution defined by Silva *et al.* (2013) and Silva (2013) |
| $$F(x) = \frac{1 - exp\left(-\lambda G(x)\right)}{1 - e^{-\lambda}}$$ | Kumaraswamy-G Poisson defined by Ramos (2014) |
| $F(x) = \left(1 - \left(1 - G^{a}(x)\right)^{b}\right)^{c}$, where $a > 0,\ b > 0$ and $c > 0$ | Kumaraswamy-G exponentiated defined by Ramos (2014) |
| $$F(x) = \frac{e^{\lambda e^{-\beta x^{\alpha}}} - e^{\lambda}}{1 - e^{\lambda}}, x > 0$$ | Beta Weibull Poisson Family defined by Paixão (2014) |
| $F(x) = \int_0^{G(x)} K t^{a-1}(1-t)^{b-1} exp(-ct)dt,$ where $a > 0,\ b > 0$ and $c \in \mathcal{R}$ | Beta Kummer generalized defined by Pescim *et al.* (2012) |
| $F(x) = \frac{e^{-\frac{\lambda}{\alpha}W(-\alpha e^{-\alpha})} - e^{-\frac{\lambda}{\alpha}W(\psi(x))}}{e^{-\frac{\lambda}{\alpha}W(-\alpha e^{-\alpha})} - 1}$, where $W(x) = \sum_{n=1}^{\infty} \frac{(-1)^{n-1} n^{n-2}}{(n-1)!} x^{n}$ and $\psi(x) = -\alpha e^{-\alpha - b x^{a}}$ | Weibull Generalized Poisson distribution defined by Paixão (2014) |
| $$F(x) = \frac{(1-\beta)^{-s} - \{1 - \beta[1 - G(x)]\}^{-s}}{(1-\beta)^{-s} - 1}$$ | G-Negative Binomial family defined by Paixão (2014) |
| $F(x) = \frac{\zeta(s) - Li_s[1 - G(x)]}{\zeta(s)}$, where $Li_s(z) = \sum_{k=1}^{\infty} \frac{z^k}{k^s}$ and $\zeta(s) = \sum_{k=1}^{\infty} \frac{1}{k^s}$ | Zeta-G defined by Paixão (2014) |
| $$F(x) = \sum_{k=0}^{x} \frac{C^{(k)}(a)}{k!\,C(\lambda)} (\lambda - a)^{k}$$ | Power Series Distributions Family (CONSUL and FAMOYE, 2006) |
| $$F(x) = \sum_{k=1}^{x} \frac{1}{k!} \left[\left(C(0)\right)^{k}\right]^{(k-1)}$$ | Basic Lagrangian defined by Consul and Famoye (2006) |
| $$F(x) = \sum_{k=n}^{x} \frac{n}{(k-n)!\,k} \left[\left(C(0)\right)^{k}\right]^{(k-n)}$$ | Lagrangian Delta defined by Consul and Famoye (2006) |
| $F(x) = (\sum_{k=0}^{x} P(X=k))^{\delta}$, $P(X=k) = \begin{cases} w(0), & k=0 \\ \left[\left(C(0)\right)^{k} w^{(1)}(0)\right]^{(k-1)}, & k=1,2,3,\ldots \end{cases}$ | Generalized Lagrangian defined by Consul and Famoye (2006) |
| $F(x) = \int_{-\infty}^{x} e^{\int \frac{a_0+a_1 t+\cdots+a_s t^s}{b_0+b_1 t+\cdots+b_r t^r} dt} dt$ | Generalized Pearson in Ordinary Differential Equation form defined by Shakil *et al.* (2010) |
| $F(x) = \int_{-\infty}^{x} e^{\int \frac{a_0+a_1 t+\cdots+a_s t^s}{b_0+b_1 t+\cdots+b_r t^r} (f(t))^{\beta} dt} dt, \quad \beta \geq 0$ | Generalized Pearson in Ordinary Differential Equation form defined by Shakil *et al.* (2010) |
| $F(x) = \int_{-\infty}^{x} \int_{-\infty}^{y} \left( \sum_{i=1}^{2} \alpha_i(t) f^{\beta_i}(t) \right) dt\, dy$ | Generalized Family in Ordinary Differential Equation form defined by (VODA, 2009) |

The aim of this work is to propose a method to create distributions and probabilistic distribution classes that could unify the various methods to generate distribution classes already described in literature. The idea of this method is to generate classes from already known distributions, using monotonic functions and a cumulative distribution function.

This paper is organized in the following way: in Section 2, we describe two methods to generate probability distributions, establishing the conditions that must be satisfied by the used monotonic functions and probability distribution to guarantee that the proposed method indeed generates a probability distribution. In Section 3, we analyze a special case of the methods described in the previous section for the case where the monotonic functions are compositions of known probability distribution functions. Still in Section 3 we present several specific cases of these methods that may be easily used to obtain new probability distributions. At the end of this section, we demonstrate that all methods presented in Sections 2 and 3 are equivalent. In Section 4, we analyze the support and nature of the distributions generated by the methods proposed in Section 3. Section 5 presents our conclusions and directions for further works. As an application of the proposed methods, the appendix to the article contains a table showing how to obtain several classes of probability distributions described in the literature using our proposed methods.

## 2. The method

The method we suggest to create distribution classes uses monotonic functions, Ʊ: $\mathbb{R}\to\mathbb{R}$, Ɣ: $\mathbb{R}\to\mathbb{R}$, $U_j$: $\mathbb{R}\to\mathbb{R}\cup\{\pm\infty\}$, $L_j$: $\mathbb{R}\to\mathbb{R}\cup\{\pm\infty\}$, $M_j$: $\mathbb{R}\to\mathbb{R}\cup\{\pm\infty\}$ and $V_j$: $\mathbb{R}\to\mathbb{R}\cup\{\pm\infty\}$, and a cumulative distribution function $F$ $(cdf)$. The idea of this method is to generate a probability distribution integrating $F$ from $L_j(x)$ to $U_j(x)$ and from $M_j(x)$ to $V_j(x)$ for any $x \in \mathbb{R}$ and $j = 1, 2, 3, \dots, n$. Theorem 1 that follows shows sufficient conditions that the functions Ʊ$(x)$, Ɣ$(x)$, $L_j(x)$, $U_j(x)$, $M_j(x)$ and $V_j(x)$ must satisfy to guarantee that the method generates a probability distribution function.

Theorem 1 (T1): Method to generate distributions and classes of probability distributions.

Let $F$: $\mathbb{R}\to\mathbb{R}$, Ʊ: $\mathbb{R}\to\mathbb{R}$, Ɣ: $\mathbb{R}\to\mathbb{R}$, $U_j$: $\mathbb{R}\to\mathbb{R}\cup\{\pm\infty\}$, $L_j$: $\mathbb{R}\to\mathbb{R}\cup\{\pm\infty\}$, $M_j$: $\mathbb{R}\to\mathbb{R}\cup\{\pm\infty\}$ and $V_j$: $\mathbb{R}\to\mathbb{R}\cup\{\pm\infty\}$, for $j = 1, 2, 3, \dots, n$, be monotonic and right continuous functions such that:

$[c1]$ $F$ is a $cdf$ and Ʊ and Ɣ are non-negative;

$[c2]$ Ʊ$(x)$, $U_j(x)$ and $M_j(x)$ are non-decreasing and Ɣ$(x)$, $V_j(x)$ and $L_j(x)$ are non-increasing $\forall\, j = 1, 2, 3, \dots, n$;

$[c3]$ If $\lim_{x\to-\infty}$ Ʊ$(x) \neq \lim_{x\to-\infty}$ Ɣ$(x)$, then $\lim_{x\to-\infty}$ Ʊ$(x) = 0$ or $\lim_{x\to-\infty} U_j(x) = \lim_{x\to-\infty} L_j(x)$ $\forall\, j = 1,2,3, \dots, n$, and $\lim_{x\to-\infty}$ Ɣ$(x) = 0$ or $\lim_{x\to-\infty} M_j(x) = \lim_{x\to-\infty} V_j(x)$, $\forall\, j = 1, 2, 3, \dots, n$;

$[c4]$ If $\lim_{x\to-\infty}$ Ʊ$(x) = \lim_{x\to-\infty}$ Ɣ$(x) \neq 0$, then $\lim_{x\to-\infty} U_j(x) = \lim_{x\to-\infty} V_j(x)$ and $\lim_{x\to-\infty} M_j(x) = \lim_{x\to-\infty} L_j(x)$, $\forall\, j = 1,2,3, \dots, n$;

$[c5]$ $\lim_{x\to-\infty} L_j(x) \leq \lim_{x\to-\infty} U_j(x)$ and if $\lim_{x\to-\infty} \mathrm{V}(x) \neq 0$, then $\lim_{x\to+\infty} M_j(x) \leq \lim_{x\to+\infty} V_j(x)$, $\forall\, j = 1,2,3,\dots,n$;

$[c6]$ $\lim_{x\to+\infty} U_n(x) \geq sup\{x \in \mathbb{R}: F(x) < 1\}$ and $\lim_{x\to+\infty} L_1(x) \leq inf\{x \in \mathbb{R} : F(x) > 0\}$;

$[c7]$ $\lim_{x\to+\infty} \mathrm{U}(x) = 1$;

$[c8]$ $\lim_{x\to+\infty} \mathrm{V}(x) = 0$ or $\lim_{x\to+\infty} M_j(x) = \lim_{x\to+\infty} V_j(x)$, $\forall\, j = 1,2,3,\dots,n$ and $n \geq 1$;

$[c9]$ $\lim_{x\to+\infty} U_j(x) = \lim_{x\to+\infty} L_{j+1}(x)$, $\forall\, j = 1,2,3,\dots,n-1$ and $n \geq 2$;

$[c10]$ $F$ is a $cdf$ without points of discontinuity or all functions $L_j(x)$ and $V_j(x)$ are constant at the right of the vicinity of points whose image are points of discontinuity of $F$, being also continuous in that points. Moreover, $F$ does not have any point of discontinuity in the set $\left\{\lim_{x\to\pm\infty} L_j(x), \lim_{x\to\pm\infty} U_j(x), \lim_{x\to\pm\infty} M_j(x), \lim_{x\to\pm\infty} V_j(x), for\ some\ j = 1,2,\dots,n\right\}$.

Then, $H(x) = \mathrm{U}(x)\sum_{j=1}^{n}\int_{L_j(x)}^{U_j(x)} dF(t) - \mathrm{V}(x)\sum_{j=1}^{n}\int_{M_j(x)}^{V_j(x)} dF(t)$ is a $cdf$.

Demonstration:

$(i)$ $\lim_{x\to-\infty} H(x) = 0$.

$$\lim_{x\to-\infty} H(x) = \lim_{x\to-\infty}\left(\mathrm{U}(x)\sum_{j=1}^{n}\int_{L_j(x)}^{U_j(x)} dF(t)\right) - \lim_{x\to-\infty}\left(\mathrm{V}(x)\sum_{j=1}^{n}\int_{M_j(x)}^{V_j(x)} dF(t)\right)$$
$$= \left(\lim_{x\to-\infty}\mathrm{U}(x)\right)\sum_{j=1}^{n}\int_{\lim_{x\to-\infty} L_j(x)}^{\lim_{x\to-\infty} U_j(x)} dF(t) - \left(\lim_{x\to-\infty}\mathrm{V}(x)\right)\sum_{j=1}^{n}\int_{\lim_{x\to-\infty} M_j(x)}^{\lim_{x\to-\infty} V_j(x)} dF(t),$$

where the last equality holds because $F$ is continuous in

$$\left\{\lim_{x\to-\infty} U_j(x), \lim_{x\to-\infty} L_j(x), \lim_{x\to-\infty} V_j(x), \lim_{x\to-\infty} M_j(x)\right\}.$$

Conditions [c3] and [c4] guarantee that:

$$\lim_{x\to-\infty} H(x) = \left(\lim_{x\to-\infty} \mho(x)\right)\sum_{j=1}^{n}\int_{\lim_{x\to-\infty} L_j(x)}^{\lim_{x\to-\infty} U_j(x)} dF(t) - \left(\lim_{x\to-\infty} \Lambda(x)\right)\sum_{j=1}^{n}\int_{\lim_{x\to-\infty} M_j(x)}^{\lim_{x\to-\infty} V_j(x)} dF(t)$$

$$= 0.$$

$(ii)$ $\lim_{x\to+\infty} H(x) = 1.$

$$\lim_{x\to+\infty} H(x) = \lim_{x\to+\infty}\left(\mho(x)\sum_{j=1}^{n}\int_{L_j(x)}^{U_j(x)} dF(t)\right) - \lim_{x\to+\infty}\left(\Lambda(x)\sum_{j=1}^{n}\int_{M_j(x)}^{V_j(x)} dF(t)\right) =$$

$$\left(\lim_{x\to+\infty} \mho(x)\right)\sum_{j=1}^{n}\int_{\lim_{x\to+\infty} L_j(x)}^{\lim_{x\to+\infty} U_j(x)} dF(t) - \left(\lim_{x\to+\infty} \Lambda(x)\right)\sum_{j=1}^{n}\int_{\lim_{x\to+\infty} M_j(x)}^{\lim_{x\to+\infty} V_j(x)} dF(t),$$

where the last equality holds because F is continuous in

$$\left\{\lim_{x\to+\infty} U_j(x), \lim_{x\to+\infty} L_j(x), \lim_{x\to+\infty} V_j(x), \lim_{x\to+\infty} M_j(x)\right\}.$$

Thus, conditions $[c1]$, $[c6]$, $[c7]$, $[c8]$ and $[c9]$ guarantee that

$$\lim_{x\to+\infty} H(x) = 1.$$

$(iii)$ If $x_1 \le x_2$, then $H(x_1) \le H(x_2)$.

Let $x_1 \le x_2$, then $[c2]$ implies that: $U_j(x_1) \le U_j(x_2)$, $L_j(x_1) \ge L_j(x_2)$, $M_j(x_1) \le M_j(x_2)$, $V_j(x_1) \ge V_j(x_2)$, $\mho(x_1) \le \mho(x_2)$ and $\Lambda(x_1) \ge \Lambda(x_2)$. Beside this, $[c2]$ and $[c5]$ imply, $\sum_{j=1}^{n}\int_{L_j(x_1)}^{U_j(x_1)} dF(t) \ge 0$, $\sum_{j=1}^{n}\int_{M_j(x_1)}^{V_j(x_1)} dF(t) \ge 0$, $\sum_{j=1}^{n}\int_{L_j(x_2)}^{U_j(x_2)} dF(t) \ge 0$ and $\sum_{j=1}^{n}\int_{M_j(x_2)}^{V_j(x_2)} dF(t) \ge 0$.

Thus, since, by $[c1]$, $\mho$ and $\Lambda$ are non-negative, we have

$$H(x_1) = \mho(x_1)\sum_{j=1}^{n}\int_{L_j(x_1)}^{U_j(x_1)} dF(t) - \Lambda(x_1)\sum_{j=1}^{n}\int_{M_j(x_1)}^{V_j(x_1)} dF(t)$$

$$\le \mho(x_2)\sum_{j=1}^{n}\int_{L_j(x_2)}^{U_j(x_2)} dF(t) - \Lambda(x_2)\sum_{j=1}^{n}\int_{M_j(x_2)}^{V_j(x_2)} dF(t) = H(x_2).$$

$(iv)$ $\lim_{x\to x_0^+} H(x) = H(x_0).$

$$\lim_{x\to x_0^+} H(x) = \lim_{x\to x_0^+} \mho(x) \sum_{j=1}^{n} \int_{L_j(x)}^{U_j(x)} dF(t) - \lim_{x\to x_0^+} \mathrm{V}(x) \sum_{j=1}^{n} \int_{M_j(x)}^{V_j(x)} dF(t)$$

$$= \left(\lim_{x\to x_0^+} \mho(x)\right) \sum_{j=1}^{n} \int_{\lim_{x\to x_0^+} L_j(x)}^{\lim_{x\to x_0^+} U_j(x)} dF(t) - \left(\lim_{x\to x_0^+} \mathrm{V}(x)\right) \sum_{j=1}^{n} \int_{\lim_{x\to x_0^+} M_j(x)}^{\lim_{x\to x_0^+} V_j(x)} dF(t)$$

$$= \mho(x_0) \sum_{j=1}^{n} \int_{L_j(x_0)}^{U_j(x_0)} dF(t) - \mathrm{V}(x_0) \sum_{j=1}^{n} \int_{M_j(x_0)}^{V_j(x_0)} dF(t) = H(x_0).$$

The above equalities hold due to $[c10]$ and because $\mho(x)$, $U_j(x)$, $M_j(x)$, $\mathrm{V}(x)$, $V_j(x)$ and $L_j(x)$ are right continuous.

From the facts $(i)$, $(ii)$, $(iii)$ and $(iv)$, we may conclude that $H(x) = \mho(x) \sum_{j=1}^{n} \int_{L_j(x)}^{U_j(x)} dF(t) - \mathrm{V}(x) \sum_{j=1}^{n} \int_{M_j(x)}^{V_j(x)} dF(t)$ is a $cdf$. ■

Corollary 1.1 presents an alternative method to generate distributions and classes of probability distributions.

Corollary 1.1 (C1.1): Complementary method to generate distributions and classes of probability distributions.

Let $\varphi: \mathbb{R}\to\mathbb{R}$, $\mho: \mathbb{R}\to\mathbb{R}$, $\mathcal{W}: \mathbb{R}\to\mathbb{R}$, $\mathbb{U}_j: \mathbb{R}\to\mathbb{R}\cup\{\pm\infty\}$, $\mathbb{L}_j: \mathbb{R}\to\mathbb{R}\cup\{\pm\infty\}$, $\mathbb{M}_j: \mathbb{R}\to\mathbb{R}\cup\{\pm\infty\}$ and $\mathbb{V}_j: \mathbb{R}\to\mathbb{R}\cup\{\pm\infty\}$, $\forall\, j = 1, 2, 3, \ldots, \eta$, be monotonic and right continuous functions such that: $[cc1]$ $F$ is a $cdf$ and $\mho$ and $\mathcal{W}$ are non-negative;

$[cc2]$ $\mho(x)$, $\mathbb{U}_j(x)$ and $\mathbb{M}_j(x)$ are non-decreasing and $\mathcal{W}(x)$, $\mathbb{V}_j(x)$ e $\mathbb{L}_j(x)$ are non-increasing, $\forall\, j = 1, 2, 3, \ldots, \eta$;

$[cc3]$ If $\lim_{x\to+\infty} \mathcal{W}(x) \neq \lim_{x\to+\infty} \mho(x)$, then $\lim_{x\to+\infty} \mho(x) = 0$ or $\lim_{x\to+\infty} \mathbb{L}_j(x) = \lim_{x\to+\infty} \mathbb{U}_j(x)$, $\forall\, j = 1,2,3,\ldots,\eta$, and $\lim_{x\to+\infty} \mathcal{W}(x) = 0$ or $\lim_{x\to+\infty} \mathbb{M}_j(x) = \lim_{x\to+\infty} \mathbb{V}_j(x)$, $\forall\, j = 1,2,3,\ldots,\eta$;

$[cc4]$ If $\lim_{x\to+\infty} \mathcal{W}(x) = \lim_{x\to+\infty} \mho(x) \neq 0$, then $\lim_{x\to+\infty} \mathbb{U}_j(x) = \lim_{x\to+\infty} \mathbb{V}_j(x)$ and $\lim_{x\to+\infty} \mathbb{M}_j(x) = \lim_{x\to+\infty} \mathbb{L}_j(x)$, $\forall\, j = 1,2,3,\ldots,\eta$;

$[cc5]$ $\lim_{x\to+\infty} \mathbb{M}_j(x) \leq \lim_{x\to+\infty} \mathbb{V}_j(x)$ and if $\lim_{x\to+\infty} \mho(x) \neq 0$, then $\lim_{x\to-\infty} \mathbb{L}_j(x) \leq \lim_{x\to-\infty} \mathbb{U}_j(x)$, $\forall\, j = 1,2,3,\ldots,\eta$;

$[cc6]$ $\lim_{x\to-\infty} \mathbb{V}_\eta(x) \geq sup\{x \in \mathbb{R}: \varphi(x) < 1\}$ and $\lim_{x\to-\infty} L_1(x) \leq inf\{x \in \mathbb{R} : \varphi(x) > 0\}$;

$[cc7]$ $\lim_{x \to -\infty} \mathcal{W}(x) = 1$;

$[cc8]$ $\lim_{x \to -\infty} \mho(x) = 0$ or $\lim_{x \to -\infty} \mathbb{L}_j(x) = \lim_{x \to -\infty} \mathbb{U}_j(x)$, $\forall\, j = 1,2,3, \dots, \eta$ and $\eta \geq 1$;

$[cc9]$ $\lim_{x \to -\infty} \mathbb{V}_j(x) = \lim_{x \to -\infty} \mathbb{M}_{j+1}(x)$, $\forall\, j = 1, 2, 3, \dots, \eta - 1$ and $\eta \geq 2$;

$[cc10]$ $\varphi$ is a $cdf$ without points of discontinuity or all functions $\mathbb{L}_j(x)$ and $\mathbb{V}_j$ are constant at the right of the vicinity of points whose image are points of discontinuity of $\varphi$, being also continuous in that points. Moreover, $\varphi$ does not have any point of discontinuity in the set $\left\{\lim_{x \to \pm\infty} \mathbb{L}_j(x)\,(x),\ \lim_{x \to \pm\infty} \mathbb{U}_j\,(x),\ \lim_{x \to \pm\infty} \mathbb{M}_j\,(x)\,,\ \lim_{x \to \pm\infty} \mathbb{V}_j\,(x),\ for\ some\ j = 1,2, \dots, \eta\right\}$.

Then, $H(x) = 1 - \mathcal{W}(x) \sum_{j=1}^{\eta} \int_{\mathbb{M}_j(x)}^{\mathbb{V}_j(x)} d\varphi(t) + \mho(x) \sum_{j=1}^{\eta} \int_{\mathbb{L}_j(x)}^{\mathbb{U}_j(x)} d\varphi(t)$ is a $cdf$.

Demonstration:

In Theorem 1, consider $n = 1$, $\mho(x) = 1$, $\mathrm{V}(x) = 0$, $U_1(x) = 1$ and $L_1(x) = \mathcal{W}(x) \sum_{j=1}^{\eta} \int_{\mathbb{M}_j(x)}^{\mathbb{V}_j(x)} d\varphi(t) - \mho(x) \sum_{j=1}^{\eta} \int_{\mathbb{L}_j(x)}^{\mathbb{U}_j(x)} d\varphi(t)$, $\forall x \in \mathbb{R}$, and $F$ a $cdf$ of the uniform [0,1] distribution. Note that $U_1(x)$ and $L_1(x)$ satisfy the hypotheses of Theorem 1, since $[cc1]$, $[cc2]$ and $[cc5]$ guarantee that $L_1(x) = \mathcal{W}(x) \sum_{j=1}^{\eta} \int_{\mathbb{M}_j(x)}^{\mathbb{V}_j(x)} d\varphi(t) - \mho(x) \sum_{j=1}^{\eta} \int_{\mathbb{L}_j(x)}^{\mathbb{U}_j(x)} d\varphi(t)$ is non-increasing and $U_1(x) = 1$ is non-decreasing. Thus, conditions $[c2]$ and $[c5]$ are satisfied. Moreover, conditions $[cc3]$ and $[cc4]$ guarantee that:

$\lim_{x \to -\infty} U_1(x) = \lim_{x \to -\infty} L_1(x) = 1$, $\lim_{x \to +\infty} U_1(x) = sup\{x \in \mathbb{R}: F(x) < 1\} = 1$, $\lim_{x \to +\infty} L_1(x) = inf\{x \in \mathbb{R} : F(x) > 0\} = 0$, that both $L_1(x)$ and $U_1(x)$ are right continuous and that $F$ is a $cdf$ without points of discontinuity.

As all conditions of Theorem 1 are satisfied, it follows that

$$H(x) = \int_{L_1(x)}^{U_1(x)} dF(s) = \int_{\mathcal{W}(x) \sum_{j=1}^{\eta} \int_{\mathbb{M}_j(x)}^{\mathbb{V}_j(x)} d\varphi(t) - \mho(x) \sum_{j=1}^{\eta} \int_{\mathbb{L}_j(x)}^{\mathbb{U}_j(x)} d\varphi(t)}^{1} ds$$

$$= 1 - \mathcal{W}(x) \sum_{j=1}^{\eta} \int_{\mathbb{M}_j(x)}^{\mathbb{V}_j(x)} d\varphi(t) + \mho(x) \sum_{j=1}^{\eta} \int_{\mathbb{L}_j(x)}^{\mathbb{U}_j(x)} d\varphi(t)$$

is a cumulative distribution function. ■

In the next section, we present some corollaries of Theorem 1 where the monotonic functions $\mathrm{U}(x)$, $\mathrm{V}(x)$, $U_j(x)$, $L_j(x)$, $M_j(x)$ and $V_j(x)$ are compositions of monotonic functions of known probability distributions.

## 3. Monotonic functions involving probabilities distributions

In this section, we show how to generate classes of probability distributions using monotonic functions which are compositions known probability distributions. Formally, consider that $\mathcal{U}: [0,1]^m \to \mathbb{R}$, $\vartheta: [0,1]^m \to \mathbb{R}$, $\mu_j: [0,1]^m \to \mathbb{R} \cup \{\pm\infty\}$, $\ell_j: [0,1]^m \to \mathbb{R} \cup \{\pm\infty\}$, $\mathcal{v}_j: [0,1]^m \to \mathbb{R}\cup\{\pm\infty\}$ and $\mathcal{m}_j: [0,1]^m \to \mathbb{R}\cup\{\pm\infty\}$ are monotonic and right continuous functions. The results of this section are achieved considering that: $\mathrm{U}(x) = \mathcal{U}(G_1, \dots, G_m)(x)$, $\mathrm{V}(x) = \vartheta(G_1, \dots, G_m)(x)$, $U_j(x) = \mu_j(G_1, \dots, G_m)(x)$, $L_j(x) = \ell_j(G_1, \dots, G_m)(x)$, $M_j(x) = \mathcal{m}_j(G_1, \dots, G_m)(x)$ and $V_j(x) = \mathcal{v}_j(G_1, \dots, G_m)(x)$.

We use the abbreviation $(.)(x) = (G_1, \dots, G_m)(x) = \big(G_1(x), \dots, G_m(x)\big)$ to represent the vector formed by the $cdf$'s calculated on the same point of the domain $x$.

Corollary 1.2 (C1.2): Method to generate classes of probability distributions.

Let $F: \mathbb{R} \to \mathbb{R}$, $\mu_j: [0,1]^m \to \mathbb{R}\cup\{\pm\infty\}$, $\ell_j: [0,1]^m \to \mathbb{R}\cup\{\pm\infty\}$, $\mathcal{U}: [0,1]^m \to \mathbb{R}$, $\mathcal{v}_j: [0,1]^m \to \mathbb{R}\cup\{\pm\infty\}$, $\mathcal{m}_j: [0,1]^m \to \mathbb{R}\cup\{\pm\infty\}$ and $\vartheta: [0,1]^m \to \mathbb{R}$, $\forall\, j = 1,2,3, \dots, n$, be monotonic and right continuous functions such that:

$[d1]$ $F$ is a $cdf$ and $\mathcal{U}$ and $\vartheta$ are non-negative;

$[d2]$ $\mu_j$, $\mathcal{m}_j$ and $\mathcal{U}$ are non-decreasing and $\ell_j$, $\mathcal{v}_j$ and $\vartheta$ are non-increasing, $\forall\, j = 1,2,3, \dots, n$, in all of its variables;

$[d3]$ If $\mathcal{U}(0, \dots, 0) \neq \vartheta(0, \dots, 0)$, then $\mathcal{U}(0, \dots, 0) = 0$ or $\mu_j(0, \dots, 0) = \ell_j(0, \dots, 0)$, $\forall\, j = 1,2,3, \dots, n$, and $\vartheta(0, \dots, 0) = 0$ or $\mathcal{m}_j(0, \dots, 0) = \mathcal{v}_j(0, \dots, 0)$, $\forall\, j = 1,2,3, \dots, n$;

$[d4]$ If $\mathcal{U}(0, \dots, 0) = \vartheta(0, \dots, 0) \neq 0$, then $\mu_j(0, \dots, 0) = \mathcal{v}_j(0, \dots, 0)$ and $\mathcal{m}_j(0, \dots, 0) = \ell_j(0, \dots, 0)$, $\forall\, j = 1,2,3, \dots, n$;

$[d5]$ $\ell_j(0,\ldots,0) \leq \mu_j(0,\ldots,0)$ and if $\vartheta(0,\ldots,0) \neq 0$, then $m_j(1,\ldots,1) \leq v_j(1,\ldots,1)$, $\forall\, j = 1,2,3,\ldots,n$;

$[d6]$ $\mu_n(1,\ldots,1) \geq sup\{x \in \mathbb{R}: F(x) < 1\}$ and $\ell_1(1,\ldots,1) \leq inf\{x \in \mathbb{R} : F(x) > 0\}$;

$[d7]$ $\mathcal{U}(1,\ldots,1) = 1$;

$[d8]$ $\vartheta(1,\ldots,1) = 0$ or $v_j(1,\ldots,1) = m_j(1,\ldots,1)$, $\forall\, j = 1,2,3,\ldots,n$ and $n \geq 1$;

$[d9]$ $\mu_j(1,\ldots,1) = \ell_{j+1}(1,\ldots,1)$, $\forall\, j = 1,2,3,\ldots,n-1$ and $n \geq 2$;

$[d10]$ $F$ is a $cdf$ without points of discontinuity or the functions $\ell_j(.)(x)$ and $v_j(.)(x)$ are constant at the right of the vicinity of points whose image are points of discontinuity of $F$, being also continuous in that points. Moreover, $F$ does not have any point of discontinuity in the set $\{\ell_j(0,\ldots,0), \mu_j(0,\ldots,0), m_j(0,\ldots,0), v_j(0,\ldots,0), \ell_j(1,\ldots,1), \mu_j(1,\ldots,1), m_j(1,\ldots,1), v_j(1,\ldots,1), for\ some\ j = 1,2,\ldots,n\}$.

Then, $H_{G_1,\ldots,G_m}(x) = \mathcal{U}(.)(x)\sum_{j=1}^{n}\int_{\ell_j(.)(x)}^{\mu_j(.)(x)} dF(t) - \vartheta(.)(x)\sum_{j=1}^{n}\int_{m_j(.)(x)}^{v_j(.)(x)} dF(t)$ is a functional generator of classes of probability distributions where $(.)(x) = (G_1,\ldots,G_m)(x)$.

Demonstration:

In Theorem 1, set $\mho(x) = \mathcal{U}(.)(x)$, $\mathrm{V}(x) = \vartheta(.)(x)$, $U_j(x) = \mu_j(.)(x)$, $L_j(x) = \ell_j(.)(x)$, $M_j(x) = m_j(.)(x)$ and $V_j(x) = v_j(.)(x)$, and observe that condition [di] implies condition [ci] of Theorem 1, for i= 1,2,….10. ■

Let us now consider a special case of Corollary 1.2 that is a functional constructor of classes of probability distributions that can be easily used.

1st special case of Corollary 1.2 (1C1.2): Easy to use method for the construction of classes of probability distributions.

Let $u_i: [0,1]^m \to [0,1]$ and $v_i: [0,1]^m \to [0,1]$ be monotonic and right continuous functions such that $u_i$'s are non-decreasing and $v_i$'s are non-increasing in each one of its variables, with $u_i(0,\ldots,0) = 0$, $u_i(1,\ldots,1) = 1$, $v_i(0,\ldots,0) = 1$ and $v_i(1,\ldots,1) = 0$, for all $i = 1,\ldots,k$. If, in Corollary 1.2, $\mathcal{U}(.)(x) = \prod_{i=1}^{k}\left((1-\theta_i)u_i(.)(x) + \theta_i\right)^{\alpha_i}$ and $\vartheta(.)(x) = \prod_{i=1}^{k}\left(\theta_i v_i(.)(x)\right)^{\alpha_i}$, with $\alpha_i \geq 0$ and $0 \leq \theta_i \leq 1$, then $H_{G_1,\ldots,G_m}(x) =$ $\prod_{i=1}^{k}\left((1-\theta_i)u_i(.)(x)+\theta_i\right)^{\alpha_i}\sum_{j=1}^{n}\int_{\ell_j(.)(x)}^{\mu_j(.)(x)} dF(t)-\prod_{i=1}^{k}\left(\theta_i v_i(.)(x)\right)^{\alpha_i}\sum_{j=1}^{n}\int_{m_j(.)(x)}^{v_j(.)(x)} dF(t)$ is a functional generator of classes of probability distributions, where $(.)(x)=(G_1,\dots,G_m)(x)$.

Table 2 shows some particular cases of the functional constructor of classes of probability distributions, $H_{G_1,\dots,G_m}(x)=\prod_{i=1}^{k}\left((1-\theta_i)u_i(.)(x)+\theta_i\right)^{\alpha_i}\sum_{j=1}^{n}\int_{\ell_j(.)(x)}^{\mu_j(.)(x)} dF(t)-\prod_{i=1}^{k}\left(\theta_i v_i(.)(x)\right)^{\alpha_i}\sum_{j=1}^{n}\int_{m_j(.)(x)}^{v_j(.)(x)} dF(t)$, that may be more easily used for the generation of classes of distribution. Consider the following functions in the expressions from **15S1C1.2** to **20S1C1.2**: $\mu:[0,1]\rightarrow\mathbb{R}\cup\{\pm\infty\}$, $\ell:[0,1]\rightarrow\mathbb{R}\cup\{\pm\infty\}$, $v:[0,1]\rightarrow\mathbb{R}\cup\{\pm\infty\}$, $m:[0,1]\rightarrow\mathbb{R}\cup\{\pm\infty\}$ such that $\mu$ and $m$ are non-decreasing and right continuous, and $v$ and $\ell$ are non-increasing and right continuous.

Table 2: Some functional constructors of classes of probability distributions obtained from **1C1.2.**

| **Some sub-cases of 1C1.2** | **Special conditions of monotonic functions and parameters** | **Functional constructor obtained** |
|---|---|---|
| **1S1C1.2** | $k=1, \alpha_1=0$ and $v_j(1,\dots,1)=m_j(1,\dots,1)$ | $H_{G_1,\dots,G_m}(x)=\sum_{j=1}^{n}\int_{\ell_j(.)(x)}^{\mu_j(.)(x)} dF(t)$ |
| **2S1C1.2** | $n=1, \theta_i=0$ and $v_j(1,\dots,1)=m_j(1,\dots,1)$ | $H_{G_1,\dots,G_m}(x)=\prod_{i=1}^{k} u_i^{\alpha_i}(.)(x)\int_{\ell_1(.)(x)}^{\mu_1(.)(x)} dF(t)$ |
| **3S1C1.2** | $n=1, k=1, \alpha_1=0$ and $v_1(1,\dots,1)=m_1(1,\dots,1)$ | $H_{G_1,\dots,G_m}(x)=\int_{\ell_1(.)(x)}^{\mu_1(.)(x)} dF(t)$ |
| **4S1C1.2** | $n=1, k=1, \alpha_1=0$, $v_j(1,\dots,1)=m_j(1,\dots,1)$ and $f(t)=\frac{1}{\mu_1(1,\dots,1)-\ell_1(1,\dots,1)}$, for $t$ in $[\ell_1(1,\dots,1),\mu_1(1,\dots,1)]$ | $H_{G_1,\dots,G_m}(x)=\frac{\mu_1(.)(x)-\ell_1(.)(x)}{\mu_1(1,\dots,1)-\ell_1(1,\dots,1)}$ |
| **5S1C1.2** | $n=1, k=1, \alpha_1=0,$<br>$\ell_1(.)(x)=\mu_1(0,\dots,0),$<br>$\mathcal{v}_1(1,\dots,1)=m_1(1,\dots,1)$ and<br>$f(t)=\frac{1}{\mu_1(1,\dots,1)-\mu_1(0,\dots,0)}$, for $t$ in<br>$[\mu_1(0,\dots,0),\mu_1(1,\dots,1)]$ | $H_{G_1,\dots,G_m}(x)=\dfrac{\mu_1(.)(x)-\mu_1(0,\dots,0)}{\mu_1(1,\dots,1)-\mu_1(0,\dots,0)}$ |
| **6S1C1.2** | $n=1, k=1, \alpha_1=0,$<br>$\mu_1(.)(x)=\ell_1(0,\dots,0),$<br>$\mathcal{v}_1(1,\dots,1)=m_1(1,\dots,1)$ and<br>$f(t)=\frac{1}{\ell_1(0,\dots,0)-\ell_1(1,\dots,1)}$, for $t$ in<br>$[\ell_1(1,\dots,1),\ell_1(0,\dots,0)]$ | $H_{G_1,\dots,G_m}(x)=\dfrac{\ell_1(.)(x)-\ell_1(0,\dots,0)}{\ell_1(1,\dots,1)-\ell_1(0,\dots,0)}$ |
| **7S1C1.2** | $k=1, \alpha_1=0$ and<br>$\sum_{j=1}^{n}\int_{\ell_j(.)(x)}^{\mu_j(.)(x)} dF(t)=1$ | $H_{G_1,\dots,G_m}(x)=1-\sum_{j=1}^{n}\int_{m_j(.)(x)}^{\mathcal{v}_j(.)(x)} dF(t)$ |
| **8S1C1.2** | $n = 1, \theta_i = 1,$<br><br>$\mu_1(.)(x) = +\infty$ and<br><br>$\ell_1(.)(x) = -\infty$ | $H_{G_1,\dots,G_m}(x) = 1 - \prod_{i=1}^{k} \left(v_i(.)(x)\right)^{\alpha_i} \int_{m_1(.)(x)}^{v_1(.)(x)} dF(t)$ |
| **9S1C1.2** | $n = 1, k = 1, \alpha_1 = 0,$<br><br>$\mu_1(.)(x) = +\infty$ and $\ell_1(.)(x) = -\infty$ | $H_{G_1,\dots,G_m}(x) = 1 - \int_{m_1(.)(x)}^{v_1(.)(x)} dF(t)$ |
| **10S1C1.2** | $n = 1, k = 1, \alpha_1 = 0,$<br><br>$\mu_1(.)(x) = +\infty, \ell_1(.)(x) = -\infty$ and<br><br>$f(t) = \frac{1}{v_1(0,\dots,0) - m_1(0,\dots,0)}$, for $t$ in<br><br>$[m_1(0,\dots,0), v_1(0,\dots,0)]$ | $H_{G_1,\dots,G_m}(x) = 1 - \frac{v_1(.)(x) - m_1(.)(x)}{v_1(0,\dots,0) - m_1(0,\dots,0)}$ |
| **11S1C1.2** | $n = 1, k = 1, \alpha_1 = 0,$<br><br>$m_1(.)(x) = v_1(1, \dots, 1),$<br><br>$\mu_1(.)(x) = +\infty, \ell_1(.)(x) = -\infty$ and<br><br>$f(t) = \frac{1}{v_1(0,\dots,0) - v_1(1,\dots,1)}$, for $t$ in<br><br>$[v_1(1, \dots, 1), v_1(0, \dots, 0)]$ | $H_{G_1,\dots,G_m}(x) = \dfrac{v_1(.)(x) - v_1(0, \dots, 0)}{v_1(1, \dots, 1) - v_1(0, \dots, 0)}$ |
| **12S1C1.2** | $n = 1, k = 1, \alpha_1 = 0,$<br><br>$v_1(.)(x) = m_1(1, \dots, 1),$<br><br>$\mu_1(.)(x) = +\infty, \ell_1(.)(x) = -\infty$ and<br><br>$f(t) = \frac{1}{m_1(1,\dots,1) - m_1(0,\dots,0)}$, for $t$ in<br><br>$[m_1(0, \dots, 0), m_1(1, \dots, 1)]$ | $H_{G_1,\dots,G_m}(x) = \dfrac{m_1(.)(x) - m_1(0, \dots, 0)}{m_1(1, \dots, 1) - m_1(0, \dots, 0)}$ |
| **13S1C1.2** | $n = 1, k = 1, \alpha_1 = 0,$ | $$H_{G_1,\dots,G_m}(x) = \int_{\ell_1(.)(x)}^{\mu_1(.)(x)} dF(t) - \int_{\mathscr{m}_1(.)(x)}^{\mathscr{v}_1(.)(x)} dF(t)$$ |
| **14S1C1.2** | $n = 1, k = 1, \alpha_1 = 0$ and<br>$f(t) = \frac{1}{\mu_1(1,\dots,1) - \ell_1(1,\dots,1) - \mathscr{v}_1(1,\dots,1) + \mathscr{m}_1(1,\dots,1)},$<br>for $t$ in $[\ell_1(1,\dots,1) + \mathscr{v}_1(1,\dots,1), \mathscr{m}_1(1,\dots,1) + \mu_1(1,\dots,1)]$. | $$H_{G_1,\dots,G_m}(x) = \frac{\mu_1(.)(x) - \ell_1(.)(x) - \mathscr{v}_1(.)(x) + \mathscr{m}_1(.)(x)}{\mu_1(1,\dots,1) - \ell_1(1,\dots,1) - \mathscr{v}_1(1,\dots,1) + \mathscr{m}_1(1,\dots,1)}$$ |
| **15S1C1.2** | $\mu_1(.)(x) = \mu\big((1-\gamma)u_{k+1}(.)(x) + \gamma\big),$<br>$\ell_1(.)(x) = \ell\big((1-\gamma)u_{k+1}(.)(x) + \gamma\big),$<br>$\mathscr{v}_1(.)(x) = \mu\big(\gamma v_{k+1}(.)(x)\big),$<br>$\mathscr{m}_1(.)(x) = \ell\big(\gamma v_{k+1}(.)(x)\big),$<br>$n = 1, \alpha_i > 0$ and $0 \le \gamma \le 1.$ | $$H_{G_1,\dots,G_m}(x) = \prod_{i=1}^{k}\big((1-\theta_i)u_i(.)(x) + \theta_i\big)^{\alpha_i} \int_{\ell_1(.)(x)}^{\mu_1(.)(x)} dF(t) - \prod_{i=n+1}^{n+k}\big(\theta_i v_i(.)(x)\big)^{\alpha_i} \int_{\mathscr{m}_1(.)(x)}^{\mathscr{v}_1(.)(x)} dF(t)$$ |
| **16S1C1.2** | $\mu_1(.)(x) = \mu\big((1-\gamma)u_{k+1}(.)(x) + \gamma\big),$<br>$\ell_1(.)(x) = -\infty,$<br>$\mathscr{v}_1(.)(x) = \mu\big(\gamma v_1(.)(x)\big),$<br>$\mathscr{m}_1(.)(x) = -\infty,$ | $$H_{G_1,\dots,G_m}(x) = \prod_{i=1}^{k}\big((1-\theta_i)u_i(.)(x) + \theta_i\big)^{\alpha_i} \int_{-\infty}^{\mu_1(.)(x)} dF(t) - \prod_{i=1}^{k}\big(\theta_i v_i(.)(x)\big)^{\alpha_i} \int_{-\infty}^{\mathscr{v}_1(.)(x)} dF(t)$$ |
| | $n = 1$ and $0 \leq \gamma \leq 1$. | |
| **17S1C1.2** | $\mu_1(.)(x) = +\infty$,<br>$\ell_1(.)(x) = \ell\big((1-\gamma)u_{k+1}(.)(x) + \gamma\big)$,<br>$\mathcal{v}_1(.)(x) = +\infty$,<br>$m_1(.)(x) = \ell\big(\gamma v_{k+1}(.)(x)\big)$,<br>$n = 1$ and $0 \leq \gamma \leq 1$. | $H_{G_1,\ldots,G_m}(x) = \prod_{i=1}^{k}\big((1-\theta_i)u_i(.)(x) + \theta_i\big)^{\alpha_i}\int_{\ell_1(.)(x)}^{+\infty} dF(t) - \prod_{i=1}^{k}\big(\theta_i v_i(.)(x)\big)^{\alpha_i}\int_{m_1(.)(x)}^{+\infty} dF(t)$ |
| **18S1C1.2** | $\mu_1(.)(x) = \mathcal{v}\big(\gamma v_{k+1}(.)(x)\big)$,<br>$\ell_1(.)(x) = m\big(\gamma v_{k+1}(.)(x)\big)$,<br>$\mathcal{v}_1(.)(x) = \mathcal{v}\big((1-\gamma)u_{k+1}(.)(x) + \gamma\big)$,<br>$m_1(.)(x) = m\big((1-\gamma)u_{k+1}(.)(x) + \gamma\big)$,<br>$n = 1$, $\alpha_i > 0$ and $0 \leq \gamma \leq 1$. | $H_{G_1,\ldots,G_m}(x) = \prod_{i=1}^{k}\big((1-\theta_i)u_i(.)(x) + \theta_i\big)^{\alpha_i}\int_{\ell_1(.)(x)}^{\mu_1(.)(x)} dF(t) - \prod_{i=n+1}^{n+k}\big(\theta_i v_i(.)(x)\big)^{\alpha_i}\int_{m_1(.)(x)}^{\mathcal{v}_1(.)(x)} dF(t)$ |
| **19S1C1.2** | $\mu_1(.)(x) = \mathscr{v}\big(\gamma v_{k+1}(.)(x)\big)$,<br>$\ell_1(.)(x) = -\infty$,<br>$\mathscr{v}_1(.)(x) = \mathscr{v}\big((1-\gamma)u_{k+1}(.)(x) + \gamma\big)$,<br>$\mathscr{m}_1(.)(x) = -\infty$,<br>$n = 1$ and $0 \le \gamma \le 1$. | $H_{G_1,\dots,G_m}(x) = \prod_{i=1}^{k}\big((1-\theta_i)u_i(.)(x) + \theta_i\big)^{\alpha_i}\int_{-\infty}^{\mu_1(.)(x)} dF(t)$<br>$- \prod_{i=1}^{k}\big(\theta_i v_i(.)(x)\big)^{\alpha_i}\int_{-\infty}^{\mathscr{v}_1(.)(x)} dF(t)$ |
| **20S1C1.2** | $\mu_1(.)(x) = +\infty$,<br>$\ell_1(.)(x) = \mathscr{m}\big(\gamma v_{k+1}(.)(x)\big)$,<br>$\mathscr{v}_1(.)(x) = +\infty$,<br>$\mathscr{m}_1(.)(x) = \mathscr{m}\big((1-\gamma)u_{k+1}(.)(x) + \gamma\big)$,<br>$n = 1$ and $0 \le \gamma \le 1$. | $H_{G_1,\dots,G_m}(x) = \prod_{i=1}^{k}\big((1-\theta_i)u_i(.)(x) + \theta_i\big)^{\alpha_i}\int_{\ell_1(.)(x)}^{+\infty} dF(t)$<br>$- \prod_{i=1}^{k}\big(\theta_i v_i(.)(x)\big)^{\alpha_i}\int_{\mathscr{m}_1(.)(x)}^{+\infty} dF(t)$ |
| **21S1C1.2** | $n = 1$. | $H_{G_1,\dots,G_m}(x) = \prod_{i=1}^{k}\big((1-\theta_i)u_i(.)(x) + \theta_i\big)^{\alpha_i}\int_{\ell_1(.)(x)}^{\mu_1(.)(x)} dF(t)$<br>$- \prod_{i=1}^{k}\big(\theta_i v_i(.)(x)\big)^{\alpha_i}\int_{\mathscr{m}_1(.)(x)}^{\mathscr{v}_1(.)(x)} dF(t)$ |
| **22S1C1.2** | $\mu_1(.)(x) = +\infty$,<br><br>$\ell_1(.)(x) = -\infty$,<br><br>$\mathcal{v}_1(.)(x) = +\infty$,<br><br>$m_1(.)(x) = -\infty$,<br><br>$n = 1$ and $\alpha_i > 0$. | $H_{G_1,\dots,G_m}(x) = \prod_{i=1}^{k}\left((1-\theta_i)u_i(.)(x) + \theta_i\right)^{\alpha_i} - \prod_{i=1}^{k}\left(\theta_i v_i(.)(x)\right)^{\alpha_i}$ |

Corollary 1.3 shows an alternative method to obtain classes of probability distributions from Corollary 1.1. It shows what hypotheses $\mathcal{U}$, $\vartheta$, $\mu_j$, $\ell_j$, $\mathcal{v}_j$ and $m_j$ must satisfy so that the functions $\mho(x)$, $\mathcal{W}(x)$, $\mathbb{U}_j(x)$, $\mathbb{L}_j(x)$, $\mathbb{M}_j(x)$ and $\mathbb{V}_j(x)$ satisfy the conditions of Corollary 1.1 and classes of probability distributions can be obtained.

Corollary 1.3 (C1.3): Complementary method to generate classes of probability distributions.

Let $\varphi: \mathbb{R} \to \mathbb{R}$, $\mu_j: [0,1]^m \to \mathbb{R} \cup \{\pm\infty\}$, $\ell_j: [0,1]^m \to \mathbb{R} \cup \{\pm\infty\}$, $\mathcal{U}: [0,1]^m \to \mathbb{R}$, $\mathcal{v}_j: [0,1]^m \to \mathbb{R} \cup \{\pm\infty\}$, $m_j: [0,1]^m \to \mathbb{R} \cup \{\pm\infty\}$ and $\vartheta: [0,1]^m \to \mathbb{R}$, $\forall\, j = 1,2,3, \dots, \eta$, be monotonic and right continuous functions such that:

$[cd1]$ $\varphi$ is a $cdf$ and $\mathcal{U}$ and $\vartheta$ are non-negative;

$[cd2]$ $\mu_j$, $m_j$ and $\mathcal{U}$ are non-decreasing and $\ell_j$, $\mathcal{v}_j$ and $\vartheta$ are non-increasing, $\forall\, j = 1,2,3, \dots, \eta$, in all of its variables;

$[cd3]$ If $\mathcal{U}(1, \dots, 1) \neq \vartheta(1, \dots, 1)$, then $\vartheta(1, \dots, 1) = 0$ or $m_j(1, \dots, 1) = \mathcal{v}_j(1, \dots, 1)$, $\forall\, j = 1, 2, 3, \dots, \eta$, and $\mathcal{U}(1, \dots, 1) = 0$ or $\ell_j(1, \dots, 1) = \mu_j(1, \dots, 1)$, $\forall\, j = 1, 2, 3, \dots, \eta$;

$[cd4]$ If $\mathcal{U}(1, \dots, 1) = \vartheta(1, \dots, 1) \neq 0$, then $\mu_j(1, \dots, 1) = \mathcal{v}_j(1, \dots, 1)$, $\forall\, j = 1,2,3, \dots, \eta$, and $m_j(1, \dots, 1) = \ell_j(1, \dots, 1)$, $\forall\, j = 1,2,3, \dots, \eta$;

$[cd5]$ $\ell_j(0, \dots, 0) \leq \mu_j(0, \dots, 0)$ and if $\vartheta(1, \dots, 1) \neq 0$, then $m_j(1, \dots, 1) \leq \mathcal{v}_j(1, \dots, 1)$, $\forall\, j = 1,2,3, \dots, \eta$;

$[cd6]$ $\mathcal{v}_\eta(0, \dots, 0) \geq sup\{x \in \mathbb{R}: F(x) < 1\}$ and $m_1(0, \dots, 0) \leq inf\{x \in \mathbb{R} : F(x) > 0\}$;

$[cd7]$ $\vartheta(0, \dots, 0) = 1$;

$[cd8]$ $\mathcal{U}(0, \dots, 0) = 0$ or $\ell_j(0, \dots, 0) = \mu_j(0, \dots, 0)$, $\forall\, j = 1, 2, 3, \dots, \eta - 1$ and $\eta \geq 1$;

$[cd9]$ $\mathcal{v}_j(0, \dots, 0) = m_{j+1}(0, \dots, 0)$, $\forall\, j = 1, 2, 3, \dots, \eta - 1$ and $\eta \geq 2$;

$[cd10]$ $\varphi$ is a $cdf$ without points of discontinuity or all functions $\ell_j(.)(x)$ and $\mathcal{v}_j(.)(x)$ are constant at the right of the vicinity of points whose image are points of discontinuity of $\varphi$, being also continuous in that points. Moreover, $\varphi$ does not have any point of discontinuity in the set $\{\ell_j(0,\dots,0), \mu_j(0,\dots,0), m_j(0,\dots,0),\ v_j(0,\dots,0), \ell_j(1,\dots,1),$ $\mu_j(1,\dots,1),\ m_j(1,\dots,1),\ v_j(1,\dots,1), para\ algum\ j = 1,2,\dots,\eta\}$.

Then, $H_{G_1,\dots,G_m}(x) = 1 - \vartheta(.)(x)\sum_{j=1}^{\eta}\int_{m_j(.)(x)}^{v_j(.)(x)} d\varphi(t) + \mathcal{U}(.)(x)\sum_{j=1}^{\eta}\int_{\ell_j(.)(x)}^{\mu_j(.)(x)} d\varphi(t)$

is a functional generator of classes of probability distributions, where $(.)(x) = (G_1,\dots,G_m)(x)$.

Demonstration:

In Corollary 1.1, set $\mho(x) = \mathcal{U}(.)(x)$, $\mathcal{W}(x) = \vartheta(.)(x)$, $\mathbb{U}_j(x) = \mu_j(.)(x)$, $\mathbb{L}_j(x) = \ell_j(.)(x)$, $\mathbb{M}_j(x) = m_j(.)(x)$ and $\mathbb{V}_j(x) = v_j(.)(x)$ and observe that condition $[cdi]$ implies condition $[cci]$ of Corollary 1.1, for $i = 1,2,\dots,10$. ■

Let us now consider a special case of Corollary 1.3 that is a functional constructor of classes of probability distributions that can be easily used.

1st special case of Corollary 1.3 (1C1.3): Easy to use complementary method for the construction of classes of probability distributions.

Let $u_i: [0,1]^m \to [0,1]$ and $v_i: [0,1]^m \to [0,1]$ be monotonic and right continuous functions such that $u_i$'s are non-decreasing and $v_i$'s are non-increasing in each one of their variables with $u_i(0,\dots,0) = 0$, $u_i(1,\dots,1) = 1$, $v_i(0,\dots,0) = 1$ and $v_i(1,\dots,1) = 0$, for all $i = 1,\dots,k$. If, in Corollary 1.3, $\vartheta(.)(x) = \prod_{i=1}^{k}\big((1-\theta_i)v_i(.)(x) + \theta_i\big)^{\alpha_i}$ and $\mathcal{U}(.)(x) = \prod_{i=1}^{k}\big(\theta_i u_i(.)(x)\big)^{\alpha_i}$, with $\alpha_i \geq 0$ and $0 \leq \theta_i \leq 1$, then $H_{G_1,\dots,G_m}(x) = 1 - \prod_{i=1}^{k}\big((1-\theta_i)v_i(.)(x) + \theta_i\big)^{\alpha_i}\sum_{j=1}^{\eta}\int_{m_j(.)(x)}^{v_j(.)(x)} d\varphi(t) + \prod_{i=1}^{k}\big(\theta_i u_i(.)(x)\big)^{\alpha_i}\sum_{j=1}^{\eta}\int_{\ell_j(.)(x)}^{\mu_j(.)(x)} d\varphi(t)$ is a functional generator of classes of probability distributions, where $(.)(x) = (G_1,\dots,G_m)(x)$.

Table 3 shows how to obtain some special cases of the function $H_{G_1,\dots,G_m}(x) = 1 - \prod_{i=1}^{k}\big((1-\theta_i)v_i(.)(x) + \theta_i\big)^{\alpha_i}\sum_{j=1}^{\eta}\int_{m_j(.)(x)}^{v_j(.)(x)} d\varphi(t) + \prod_{i=1}^{k}\big(\theta_i u_i(.)(x)\big)^{\alpha_i}\sum_{j=1}^{\eta}\int_{\ell_j(.)(x)}^{\mu_j(.)(x)} d\varphi(t)$ that may be more easily used to generate classes of distributions. It is important to emphasize that we can obtain the same constructors from 1S1C1.2 to 12S1C1.2 using 1C1.3, we omit the details here showing only how to obtain different constructors from those of Table 2. Consider the following functions in the expressions from **15S1C1.3** to **20S1C1.3**: $\mu: [0,1] \rightarrow \mathbb{R} \cup \{\pm\infty\}$, $\ell: [0,1] \rightarrow \mathbb{R} \cup \{\pm\infty\}$, $\mathcal{v}: [0,1] \rightarrow \mathbb{R} \cup \{\pm\infty\}$, $\mathcal{m}: [0,1] \rightarrow \mathbb{R} \cup \{\pm\infty\}$ such that $\mu$ and $\mathcal{m}$ are non-decreasing and right continuous, and $\mathcal{v}$ and $\ell$ are non-increasing and right continuous.

Table 3: Some functional constructors of classes of probability distributions obtained from **1C1.3.**

| **Some sub cases of 1C1.3** | **Special conditions over monotonic functions and parameters** | **Constructor functions obtained** |
|---|---|---|
| **13S1C1.3** | $\eta = 1, k = 1, \alpha_1 = 0,$ | $$H_{G_1,\ldots,G_m}(x) = 1 - \int_{m_1(.)(x)}^{v_1(.)(x)} dF(t) + \int_{\ell_1(.)(x)}^{\mu_1(.)(x)} dF(t)$$ |
| **14S1C1.3** | $\eta = 1, k = 1, \alpha_1 = 0$ and $f(t) = \frac{1}{v_1(0,\ldots,0) - m_1(0,\ldots,0) - \mu_1(0,\ldots,0) + \ell_1(0,\ldots,0)}$, for $t$ in $[m_1(0,\ldots,0) + \mu_1(0,\ldots,0), \ell_1(0,\ldots,0) + v_1(0,\ldots,0)]$ | $$H_{G_1,\ldots,G_m}(x) = 1 - \frac{v_1(.)(x) - m_1(.)(x) - \mu_1(.)(x) + \ell_1(.)(x)}{v_1(0,\ldots,0) - m_1(0,\ldots,0) - \mu_1(0,\ldots,0) + \ell_1(0,\ldots,0)}$$ |
| **15S1C1.3** | $\mu_1(.)(x) = \mu(\gamma u_{k+1}(.)(x))$,<br>$\ell_1(.)(x) = \ell(\gamma u_{k+1}(.)(x))$,<br>$v_1(.)(x) = \mu((1-\gamma)v_{k+1}(.)(x) + \gamma)$,<br>$m_1(.)(x) = \ell((1-\gamma)v_{k+1}(.)(x) + \gamma)$,<br>$\eta = 1$ and $0 \leq \gamma \leq 1$. | $$H_{G_1,\ldots,G_m}(x) = 1 - \prod_{i=1}^{k} \left((1-\theta_i)v_i(.)(x) + \theta_i\right)^{\alpha_i} \int_{m_1(.)(x)}^{v_1(.)(x)} d\varphi(t) + \prod_{i=1}^{k} \left(\theta_i u_i(.)(x)\right)^{\alpha_i} \int_{\ell_1(.)(x)}^{\mu_1(.)(x)} d\varphi(t)$$ |
| **16S1C1.3** | $\mu_1(.)(x) = \mu(\gamma u_{k+1}(.)(x))$,<br>$\ell_1(.)(x) = -\infty$,<br>$v_1(.)(x) = \mu((1-\gamma)v_{k+1}(.)(x) + \gamma)$,<br>$m_1(.)(x) = -\infty$,<br>$\eta = 1$ and $0 \leq \gamma \leq 1$. | $$H_{G_1,\ldots,G_m}(x) = 1 - \prod_{i=1}^{k} \left((1-\theta_i)v_i(.)(x) + \theta_i\right)^{\alpha_i} \int_{-\infty}^{v_1(.)(x)} d\varphi(t) + \prod_{i=1}^{k} \left(\theta_i u_i(.)(x)\right)^{\alpha_i} \int_{-\infty}^{\mu_1(.)(x)} d\varphi(t)$$ |
| **17S1C1.3** | $\mu_1(.)(x) = +\infty$,<br>$\ell_1(.)(x) = \ell(\gamma u_{k+1}(.)(x))$,<br>$v_1(.)(x) = +\infty$,<br>$m_1(.)(x) = \ell((1-\gamma)v_{k+1}(.)(x) + \gamma)$,<br>$\eta = 1$ and $0 \leq \gamma \leq 1$. | $H_{G_1,\dots,G_m}(x) = 1 - \prod_{i=1}^{k}((1-\theta_i)v_i(.)(x) + \theta_i)^{\alpha_i} \int_{m_1(.)(x)}^{+\infty} d\varphi(t) + \prod_{i=1}^{k}(\theta_i u_i(.)(x))^{\alpha_i} \int_{\ell_1(.)(x)}^{+\infty} d\varphi(t)$ |
| **18S1C1.3** | $\mu_1(.)(x) = m(\gamma u_{k+1}(.)(x))$,<br>$\ell_1(.)(x) = v(\gamma u_{k+1}(.)(x))$,<br>$v_1(.)(x) = m((1-\gamma)v_{k+1}(.)(x) + \gamma)$,<br>$m_1(.)(x) = v((1-\gamma)v_{k+1}(.)(x) + \gamma)$,<br>$\eta = 1$ and $0 \leq \gamma \leq 1$. | $H_{G_1,\dots,G_m}(x) = 1 - \prod_{i=1}^{k}((1-\theta_i)v_i(.)(x) + \theta_i)^{\alpha_i} \int_{m_1(.)(x)}^{v_1(.)(x)} d\varphi(t) + \prod_{i=1}^{k}(\theta_i u_i(.)(x))^{\alpha_i} \int_{\ell_1(.)(x)}^{\mu_1(.)(x)} d\varphi(t)$ |
| **19S1C1.3** | $\mu_1(.)(x) = m(\gamma u_{k+1}(.)(x))$,<br>$\ell_1(.)(x) = -\infty$,<br>$v_1(.)(x) = m((1-\gamma)v_{k+1}(.)(x) + \gamma)$,<br>$m_1(.)(x) = -\infty$,<br>$\eta = 1$ and $0 \leq \gamma \leq 1$. | $H_{G_1,\dots,G_m}(x) = 1 - \prod_{i=1}^{k}((1-\theta_i)v_i(.)(x) + \theta_i)^{\alpha_i} \int_{-\infty}^{v_1(.)(x)} d\varphi(t) + \prod_{i=1}^{k}(\theta_i u_i(.)(x))^{\alpha_i} \int_{-\infty}^{\mu_1(.)(x)} d\varphi(t)$ |
| **20S1C1.3** | $\mu_1(.)(x) = +\infty$,<br>$\ell_1(.)(x) = v(\gamma u_{k+1}(.)(x))$,<br>$v_1(.)(x) = +\infty$,<br>$m_1(.)(x) = v((1-\gamma)v_{k+1}(.)(x) + \gamma)$,<br>$\eta = 1$ and $0 \leq \gamma \leq 1$. | $H_{G_1,\dots,G_m}(x) = 1 - \prod_{i=1}^{k}((1-\theta_i)v_i(.)(x) + \theta_i)^{\alpha_i} \int_{m_1(.)(x)}^{+\infty} d\varphi(t) + \prod_{i=1}^{k}(\theta_i u_i(.)(x))^{\alpha_i} \int_{\ell_1(.)(x)}^{+\infty} d\varphi(t)$ |
| **21S1C1.3** | $\eta = 1.$ | $$H_{G_1,\ldots,G_m}(x) = 1 - \prod_{i=1}^{k}\big((1-\theta_i)v_i(.)(x)+\theta_i\big)^{\alpha_i}\int_{m_1(.)(x)}^{v_1(.)(x)} d\varphi(t) + \prod_{i=1}^{k}\big(\theta_i u_i(.)(x)\big)^{\alpha_i}\int_{\ell_1(.)(x)}^{\mu_1(.)(x)} d\varphi(t)$$ |
| **22S1C1.3** | $\mu_1(.)(x) = +\infty$,<br>$\ell_1(.)(x) = -\infty$,<br>$v_1(.)(x) = +\infty$,<br>$m_1(.)(x) = -\infty$,<br>$\eta = 1$ and $\alpha_i > 0$. | $$H_{G_1,\ldots,G_m}(x) = 1 - \prod_{i=1}^{k}\big((1-\theta_i)v_i(.)(x)+\theta_i\big)^{\alpha_i} + \prod_{i=1}^{k}\big(\theta_i u_i(.)(x)\big)^{\alpha_i}$$ |

The following theorem shows that Theorem 1 and of its corollaries are equivalent. In other words, Theorem 1 and all of its corollaries generate the same probabilistic distributions.

Theorem 2 (T2): Equivalence among Theorem 1 and its corollaries.

Theorem 1 and all of its corollaries generate exactly the same probabilistic distributions

Demonstration:

To demonstrate Theorem 2 we show that C1.1 is a corollary of T1, that C1.2 is a corollary of C1.1, that C1.3 is a corollary of C1.2, and finally that T1 is a corollary of C1.3.

(1) C1.1 is a corollary of T1: it is obvious, as it has been already demonstrated.
(2) C1.2 is a corollary of C1.1: In Corollary 1.1, $\mathcal{W}(x) = 1$, $\eta = 1$, $\mathbb{V}_1(x) = 1$, $\mathbb{M}_1(x) = \mathcal{U}(.)(x) \sum_{j=1}^{n} \int_{\ell_j(.)(x)}^{\mu_j(.)(x)} dF(t) - \vartheta(.)(x) \sum_{j=1}^{n} \int_{m_j(.)(x)}^{v_j(.)(x)} dF(t)$, $\mho(x) = 0$, $\forall x \in \mathbb{R}$ and $\varphi(t)$ the $cdf$ of the uniform[0,1].
(3) C1.3 is a corollary of C1.2: In Corollary 1.2, set $n = 1$, $\mathcal{U}(.)(x) = 1$, $\mu_1(.)(x) = 1$, $\ell_1(.)(x) = 0$, $\vartheta(.)(x) = 1$, $m_1(.)(x) = 0$, $v_1(.)(x) = \vartheta(.)(x) \sum_{j=1}^{\eta} \int_{m_j(.)(x)}^{v_j(.)(x)} d\varphi(t) - \mathcal{U}(.)(x) \sum_{j=1}^{\eta} \int_{\ell_j(.)(x)}^{\mu_j(.)(x)} d\varphi(t)$, $\forall x \in \mathbb{R}$ and $F(t)$ as the $cdf$ of the uniform [0,1].
(4) T1 is a corollary of C1.3: In Corollary 1.3, set $\eta = 1$, $\mathcal{U}(.)(x) = 1$, $\vartheta(.)(x) = 1$, $v_1(.)(x) = 1$, $m_1(.)(x) = 0$, $\ell_j(.)(x) = 0$, $\mu_1(.)(x) = \mathrm{U}(x) \sum_{j=1}^{n} \int_{L_j(x)}^{U_j(x)} dF(t) - \mathrm{V}(x) \sum_{j=1}^{n} \int_{M_j(x)}^{V_j(x)} dF(t)$, $\forall x \in \mathbb{R}$, and $\varphi(t)$ $fda$ of the uniform [0,1]. ■

From (1) to (4), we conclude Theorem 2. ■

Several classes of probability distributions existing in the literature can be obtained as special cases of the functional constructors of classes of probability distributions proposed here. Table 4, in the Appendix, shows how to obtain such classes using some corollaries of Theorem 1.

## 4. Support of the classes of probability distributions

In this section, we provide an analysis about the support and nature of the probability distributions generated through the methods described in Corollaries 1.2 and 1.3. These results are important to gain a deeper understanding about the proposed method, especially considering the fact that there is little work on this theme in the literature.

In order to state the results, we remind the reader that, by the definition, the support of a cumulative distribution function $F$ is given by $S_F = \{x \in \mathbb{R}: F(x) - F(x-\varepsilon) > 0, \forall \varepsilon > 0\}$. Theorem 3 shows that the support of the generated distribution is contained in the union of the supports of the baseline distributions $G_i$'s.

Theorem 3 (T3): General theorem of the supports.

Let $H_{G_1,\ldots,G_m}(x)$ be a cumulative distribution function generated from Corollary 1.2 (respectively, 1.3). Then, $S_{H_{G_1,\ldots,G_m}} \subset \bigcup_{j=1}^{m} S_{G_j}$.

Demonstration:

Consider that $H_{G_1,\ldots,G_m}(x)$ has the functional form of Corollary 1.2:

$$H_{G_1,\ldots,G_m}(x) = \mathcal{U}(.)(x)\sum_{j=1}^{n}\int_{\ell_j(.)(x)}^{\mu_j(.)(x)} dF(t) - \vartheta(.)(x)\sum_{j=1}^{n}\int_{m_j(.)(x)}^{\upsilon_j(.)(x)} dF(t).$$

Thus, it follows that

$$H_{G_1,\ldots,G_m}(x) - H_{G_1,\ldots,G_m}(x-\varepsilon) = \mathcal{U}(.)(x)\sum_{j=1}^{n}\int_{\ell_j(.)(x)}^{\mu_j(.)(x)} dF(t) - \vartheta(.)(x)\sum_{j=1}^{n}\int_{m_j(.)(x)}^{\upsilon_j(.)(x)} dF(t)$$

$$-\mathcal{U}(.)(x-\varepsilon)\sum_{j=1}^{n}\int_{\ell_j(.)(x-\varepsilon)}^{\mu_j(.)(x-\varepsilon)} dF(t) + \vartheta(.)(x-\varepsilon)\sum_{j=1}^{n}\int_{m_j(.)(x-\varepsilon)}^{\upsilon_j(.)(x-\varepsilon)} dF(t).$$

Suppose that $x \notin \bigcup_{j=1}^{m} S_{G_j}$. Then, there exists $\varepsilon > 0$ such that $G_j(x) - G_j(x-\varepsilon) = 0$, for all $j = 1, 2, \ldots, m$. Let us show that $H_{G_1,\ldots,G_m}(x) - H_{G_1,\ldots,G_m}(x-\varepsilon) = 0$. First, note that

$$H_{G_1,\ldots,G_m}(x) - H_{G_1,\ldots,G_m}(x-\varepsilon) = \big(\mathcal{U}(.)(x) - \mathcal{U}(.)(x-\varepsilon)\big)\sum_{j=1}^{n}\int_{\ell_j(.)(x)}^{\mu_j(.)(x)} dF(t)$$

$$-\big(\vartheta(.)(x) - \vartheta(.)(x-\varepsilon)\big)\sum_{j=1}^{n}\int_{m_j(.)(x)}^{\upsilon_j(.)(x)} dF(t)$$

$$+\mathcal{U}(.)(x-\varepsilon)\left(\sum_{j=1}^{n}\left(\int_{\mu_j(.)(x-\varepsilon)}^{\mu_j(.)(x)} dF(t) - \int_{\ell_j(.)(x-\varepsilon)}^{\ell_j(.)(x)} dF(t)\right)\right)$$

$$-\vartheta(.)(x-\varepsilon)\left(\sum_{j=1}^{n}\left(\int_{\mathcal{v}_j(.)(x-\varepsilon)}^{\mathcal{v}_j(.)(x)} dF(t) - \int_{\mathcal{m}_j(.)(x-\varepsilon)}^{\mathcal{m}_j(.)(x)} dF(t)\right)\right).$$

Since $\mathcal{U}(.)(x) = \mathcal{U}(.)(x-\varepsilon)$, $\vartheta(.)(x) = \vartheta(.)(x-\varepsilon)$, $\mu_j(.)(x) = \mu_j(.)(x-\varepsilon)$, $\ell_j(.)(x) = \ell_j(.)(x-\varepsilon)$, $\mathcal{m}_j(.)(x) = \mathcal{m}_j(.)(x-\varepsilon)$ and $\mathcal{v}_j(.)(x) = \mathcal{v}_j(.)(x-\varepsilon)$, it follows that

$$H_{G_1,\ldots,G_m}(x) - H_{G_1,\ldots,G_m}(x-\varepsilon) = 0.$$

Thus, we have $x \notin S_{H_{G_1,\ldots,G_m}}$. Therefore, $S_{H_{G_1,\ldots,G_m}} \subset \bigcup_{j=1}^{m} S_{G_j}$ . A similar argument works for the case where $H_{G_1,\ldots,G_m}(x)$ has the functional form of Corollary 1.3.

■

Corollary 3.1 (C3.1) shows a special case where the distribution $H_{G_1,\ldots,G_m}(x)$ is discrete.

Corollary 3.1 (C3.1): Discrete baselines generate discrete distributions.

If all $G_j$'s are discrete in Corollary 1.2 (respectively, 1.3), then $H_{G_1,\ldots,G_m}(x)$ is discrete.

Demonstration:

Being all $G_j$'s discrete, then $\bigcup_{j=1}^{m} S_{G_j}$ has a countable number of values. Since, by Theorem 3, $S_{H_{G_1,\ldots,G_m}} \subset \bigcup_{j=1}^{m} S_{G_j}$, then $S_{H_{G_1,\ldots,G_m}}$ has a countable number of values and, for this reason, $H_{G_1,\ldots,G_m}(x)$ is a$cdf$ of a discrete random variable. ■

Theorem 4 shows conditions when $S_{H_{G_1,\ldots,G_m}} = \bigcup_{j=1}^{m} S_{G_j}$

Theorem 4 (T4): The support of distribution is a union of the supports of the baselines.

Assume, in Corollary 1.2, (respectively, 1.3) that:

$[f1]$ $S_F$ is a convex set;

$[f2]$ $\mu_n(1,\dots,1) = sup\{x \in \mathbb{R} : \varphi(x) < 1\}$, $\ell_1(1,\dots,1) = inf\{x \in \mathbb{R} : \varphi(x) > 0\}$, $\mathcal{U}(.)(x) > 0$, $\forall\, x \in \mathbb{R}$, and that $\mu_j(.)(x)$ or $\ell_j(.)(x)$, for some $j = 1,2,\dots,n$, are strictly monotonic or $\mathcal{v}_n(0,\dots,0) = sup\{x \in \mathbb{R} : \varphi(x) < 1\}$, $\mathcal{m}_1(0,\dots,0) = inf\{x \in \mathbb{R} : \varphi(x) > 0\}$, $\vartheta(.)(x) > 0$, $\forall\, x \in \mathbb{R}$, and that $\mathcal{v}_j(.)(x)$ or $\mathcal{m}_j(.)(x)$, for some $j = 1,2,\dots,n$, are strictly monotonic.

Then, $S_{H_{G_1,\dots,G_m}} = \bigcup_{j=1}^{m} S_{G_j}$.

Demonstration:

As a consequence of Theorem 3, we only need to show that $S_{H_{G_1,\dots,G_m}} \supset \bigcup_{j=1}^{m} S_{G_j}$. Consider that $H_{G_1,\dots,G_m}(x)$ has the functional form of Corollary 1.2:

$$H_{G_1,\dots,G_m}(x) = \mathcal{U}(.)(x)\sum_{j=1}^{n}\int_{\ell_j(.)(x)}^{\mu_j(.)(x)} dF(t) - \vartheta(.)(x)\sum_{j=1}^{n}\int_{\mathcal{m}_j(.)(x)}^{\mathcal{v}_j(.)(x)} dF(t).$$

Thus, using an argument identical to that of the proof of Theorem 3, we have

$$H_{G_1,\dots,G_m}(x) - H_{G_1,\dots,G_m}(x-\varepsilon) = \big(\mathcal{U}(.)(x) - \mathcal{U}(.)(x-\varepsilon)\big)\sum_{j=1}^{n}\int_{\ell_j(.)(x)}^{\mu_j(.)(x)} dF(t)$$

$$-\big(\vartheta(.)(x) - \vartheta(.)(x-\varepsilon)\big)\sum_{j=1}^{n}\int_{\mathcal{m}_j(.)(x)}^{\mathcal{v}_j(.)(x)} dF(t)$$

$$+\mathcal{U}(.)(x-\varepsilon)\left(\sum_{j=1}^{n}\left(\int_{\mu_j(.)(x-\varepsilon)}^{\mu_j(.)(x)} dF(t) - \int_{\ell_j(.)(x-\varepsilon)}^{\ell_j(.)(x)} dF(t)\right)\right)$$

$$-\vartheta(.)(x-\varepsilon)\left(\sum_{j=1}^{n}\left(\int_{\mathcal{v}_j(.)(x-\varepsilon)}^{\mathcal{v}_j(.)(x)} dF(t) - \int_{\mathcal{m}_j(.)(x-\varepsilon)}^{\mathcal{m}_j(.)(x)} dF(t)\right)\right).$$

Suppose that $x \in \bigcup_{j=1}^{m} S_{G_j}$. Then, there exists $\varepsilon > 0$ such that $G_j(x) - G_j(x-\varepsilon) > 0$ for some $j = 1,2,\dots,m$.. Conditions $[f1]$ and $[f2]$ imply that at least one of the integrals of the form $\int_{h(.)(x-\varepsilon)}^{h(.)(x)} dF(t)$ is different from zero for $h = \mu_j$, $h = \ell_j$, $h = \mathcal{v}_j$ or $h = \mathcal{m}_j$, for some $j = 1,2,\dots,n$. This fact together with the fact that $\mathcal{U}$ or $\vartheta$ are strictly monotonic imply that:

$H_{G_1,...,G_m}(x) - H_{G_1,...,G_m}(x-\varepsilon) > 0$. Thus, $x \in S_{H_{G_1,...,G_m}}$, as desired. A similar argument works for the case where $H_{G_1,...,G_m}(x)$ has the functional form of Corollary 1.3. ■

Theorem 5 shows some conditions that guarantee that $H_{G_1,...,G_m}(x)$ is a continuous $cdf$.

Theorem 5 (T5): Generating continuous cumulative distribution functions.

Suppose that $F(x)$, $G_1$,...,$G_m$ are continuous $cdf$'s and that $\mu_j$, $\ell_j$, $\mathcal{U}$, $\mathcal{v}_j$, $\mathcal{m}_j$ and $\vartheta$ are continuous functions in Corollary 1.2 (respectively, 1.3). Then, $H_{G_1,...,G_m}(x)$ is a continuous $cdf$.

Demonstration:

Consider that $H_{G_1,...,G_m}(x)$ has the functional form of Corollary 1.2

$$H_{G_1,...,G_m}(x) = \mathcal{U}(.)(x)\sum_{j=1}^{n}\int_{\ell_j(.)(x)}^{\mu_j(.)(x)} dF(t) - \vartheta(.)(x)\sum_{j=1}^{n}\int_{\mathcal{m}_j(.)(x)}^{\mathcal{v}_j(.)(x)} dF(t).$$

Thus, it follows that

$$H_{G_1,...,G_m}(x) - H_{G_1,...,G_m}(x^-) = \mathcal{U}(.)(x)\sum_{j=1}^{n}\int_{\ell_j(.)(x)}^{\mu_j(.)(x)} dF(t) - \vartheta(.)(x)\sum_{j=1}^{n}\int_{\mathcal{m}_j(.)(x)}^{\mathcal{v}_j(.)(x)} dF(t)$$

$$-\mathcal{U}(.)(x^-)\sum_{j=1}^{n}\int_{\ell_j(.)(x^-)}^{\mu_j(.)(x^-)} dF(t) + \vartheta(.)(x^-)\sum_{j=1}^{n}\int_{\mathcal{m}_j(.)(x^-)}^{\mathcal{v}_j(.)(x^-)} dF(t).$$

Using a similar approach as that used in the proof of Theorem 3, we obtain:

$$H_{G_1,...,G_m}(x) - H_{G_1,...,G_m}(x^-) = \left(\mathcal{U}(.)(x) - \mathcal{U}(.)(x^-)\right)\sum_{j=1}^{n}\int_{\ell_j(.)(x)}^{\mu_j(.)(x)} dF(t)$$

$$-\left(\vartheta(.)(x) - \vartheta(.)(x^-)\right)\sum_{j=1}^{n}\int_{\mathcal{m}_j(.)(x)}^{\mathcal{v}_j(.)(x)} dF(t)$$

$$+\mathcal{U}(.)(x^-)\left(\sum_{j=1}^{n}\left(\int_{\mu_j(.)(x^-)}^{\mu_j(.)(x)} dF(t) - \int_{\ell_j(.)(x^-)}^{\ell_j(.)(x)} dF(t)\right)\right)$$

$$-\vartheta(.)(x^-)\left(\sum_{j=1}^{n}\int_{\mathcal{v}_j(.)(x^-)}^{\mathcal{v}_j(.)(x)} dF(t) - \sum_{j=1}^{n}\int_{m_j(.)(x^-)}^{m_j(.)(x)} dF(t)\right).$$

Since all functions included in the previous expression are continuous, we have:

$$H_{G_1,\dots,G_m}(x) - H_{G_1,\dots,G_m}(x^-) = 0.$$

Therefore, we shall conclude that $H_{G_1,\dots,G_m}(x)$ is a continuous function. A similar argument works for the case where $H_{G_1,\dots,G_m}(x)$ has the functional form of Corollary 1.3.

■

Theorem 6 shows conditions where distribution $H_{G_1,\dots,G_m}(x)$ will be a continuous $fda$ of $v.a.$.

Theorem 6 (T6): Generating cumulative distribution functions of continuous random variables.

Suppose that $F(x)$, $G_1$,..., $G_m$ are $cdf$'s of continuous random variables and that $\mu_j$, $\ell_j$, $\mathcal{U}$, $\mathcal{v}_j$, $m_j$ and $\vartheta$ are continuous and differentiable functions in Corollary 1.2 (respectively, 1.3). Then, $H_{G_1,\dots,G_m}(x)$ is a $cdf$ of a continuous random variable.

Demonstration:

Consider that $H_{G_1,\dots,G_m}(x)$ has the functional form of Corollary 1.2

$$H_{G_1,\dots,G_m}(x) = \mathcal{U}(.)(x)\sum_{j=1}^{n}\int_{\ell_j(.)(x)}^{\mu_j(.)(x)} dF(t) - \vartheta(.)(x)\sum_{j=1}^{n}\int_{m_j(.)(x)}^{\mathcal{v}_j(.)(x)} dF(t).$$

Since $F(x)$, $G_1$,..., $G_m$ are $cdf$'s of continuous random variables and $\mu_j$, $\ell_j$, $\mathcal{U}$, $\mathcal{v}_j$, $m_j$ and $\vartheta$ are continuous and differentiable functions , then $H_{G_1,\dots,G_m}(x)$ is a $cdf$ of a continuous random variable with density given by:

$$h_{G_1,\dots,G_m}(x) = \left(\sum_{z=1}^{m}\frac{\partial \mathcal{U}(.)(x)}{\partial G_z} g_z(x)\right)\left(\sum_{j=1}^{n}\int_{\ell_j(.)(x)}^{\mu_j(.)(x)} dF(t)\right)$$

$$+\mathcal{U}(.)(x)\left(\sum_{j=1}^{n}\left(F\left(\mu_j(.)(x)\right)\sum_{z=1}^{m}\frac{\partial \mu_j(.)(x)}{\partial G_z} g_z(x) - F\left(\ell_j(.)(x)\right)\sum_{z=1}^{m}\frac{\partial \ell_j(.)(x)}{\partial G_z} g_z(x)\right)\right)$$

$$-\left(\sum_{z=1}^{m}\frac{\partial\vartheta(.)(x)}{\partial G_z}g_z(x)\right)\left(\sum_{j=1}^{n}\int_{m_j(.)(x)}^{v_j(.)(x)}dF(t)\right)$$

$$-\vartheta(.)(x)\left(\sum_{j=1}^{n}\left(F\left(v_j(.)(x)\right)\sum_{z=1}^{m}\frac{\partial v_j(.)(x)}{\partial G_z}g_z(x)-F\left(m_j(.)(x)\right)\sum_{z=1}^{m}\frac{\partial m_j(.)(x)}{\partial G_z}g_z(x)\right)\right),$$

where $(.)(x)=(G_1,\dots,G_m)(x)$. A similar argument works for the case where $H_{G_1,\dots,G_m}(x)$ has the functional form of Corollary 1.3. ■

The next theorem shows an alternative way of generating discrete distributions.

Theorem 7 (T7): Integrals with respect to discrete distributions generate discrete distributions.

Suppose that the probability distribution $F(x)$ is discrete and that $\mathcal{U}(.)(x)=\vartheta(.)(x)=1$, in Corollary C1.2 (respectively, 1.3). Then, $H_{G_1,\dots,G_m}(x)$ is discrete independent from the monotonic functions that are used as limits of the integration.

Demonstration:

$$H_{G_1,\dots,G_m}(x)=\sum_{j=1}^{n}\int_{\ell_j(.)(x)}^{\mu_j(.)(x)}dF(t)-\sum_{j=1}^{n}\int_{m_j(.)(x)}^{v_j(.)(x)}dF(t)$$

$$H_{G_1,\dots,G_m}(x)=\sum_{j=1}^{n}\left(F\left(\mu_j(.)(x)\right)-F\left(\ell_j(.)(x)\right)\right)-\sum_{j=1}^{n}\left(F\left(v_j(.)(x)\right)-F\left(m_j(.)(x)\right)\right)$$

Since $F(x)$is a $cdf$ of a discrete random variable, then $F(x)$ assumes a countable number of different values. Thus, as $H_{G_1,\dots,G_m}(x)$ is given by the sum of differences of $F(x)$ evaluated in at most *4n* distinct points, it also can only assume a countable number of different values and, therefore, it is a $cdf$ of a discrete random variable.

## 5. Conclusions

The method to generate distributions and classes of probability distributions that we presented in this paper combines several methods to generate classes of distribution that have already been described in the literature. By this unification, we could draw conclusions on the supports of generated classes. Using the proposed method, we can generate any probability distribution in different ways. The only necessity is to modify the monotonic functions involved in the method.

As a further step, we aim to explore several classes of distribution which may be generated using this method, developing their properties and applying them to model several datasets. In a parallel work, we are proposing a method for generating multivariate distributions.

It is important to stress that a model can better describe a phenomenon by increasing its number of parameters, providing higher flexibility. On the other hand, we should not forget that increasing the number of parameters may cause identifiability and computational problems in the estimation of the parameters. Moreover, a large number of parameters increase the chance of overfitting, which is a problem particularly in forecasting and prediction studies. Thus, the best approach is to choose a method that best describes the analyzed phenomenon or experiment with the lowest possible number of parameters.

## Appendix: how to obtain generalizations of already existing class models

In this appendix we will show some applications to obtain some very special examples of functions generating classes of probabilistic distributions by finding probability distribution classes that have already been described in literature.

Table 4 shows how to obtain classes of probability distributions already existing in the literature by the use of some corollaries of Theorem 1

Table 4: Generalizations of already existing classes of probability distributions

| **Sub-case of 1C1.2 used** | **Used distributions** $f(t)$ | **Monotonic functions** | **Some special values for the parameters** | **Obtained class** |
|---|---|---|---|---|
| **3S1C1.2** | $\frac{1}{B(a,b)} t^{a-1}(1-t)^{b-1}$ | $\ell_1(.)(x) = \theta \prod_{i=1}^{m} \left(1 - G_i^{\alpha_i}(x)\right)^{\delta_i}$ <br> $\mu_1(.)(x) = (1-\theta) \prod_{j=1}^{m} G_j^{\beta_j}(x) + \theta$ | $\theta = 0,\ m = 1$ and $\beta_1 = 1$ | beta1 generalized defined by Eugene *et al.* (2002) |
| | | | $\theta = 0$ and $m = 1$ | Mc1 generalized defined by McDonald (1984) |
| **9S1C1.2** | $\frac{1}{B(a,b)} t^{a-1}(1-t)^{b-1}$ | $m_1(.)(x) = \theta \prod_{j=1}^{m} G_j^{\beta_j}(x)$ <br> $v_1(.)(x) = (1-\theta) \prod_{i=1}^{m} \left(1 - G_i^{\alpha_i}(x)\right)^{\delta_i} + \theta$ | $\theta = 1,\ m = 1$ and $\beta_1 = 1$ | beta1 generalized defined by Eugene *et al.* (2002) |
| | | | $\theta = 1$ and $m = 1$ | Mc1 generalized defined by McDonald (1984) |
| **3S1C1.2** | $bt^{b-1}$ | $\ell_1(.)(x) = \theta \prod_{i=1}^{m} \left(1 - G_i^{\alpha_i}(x)\right)^{\delta_i}$ | $\theta = 0,\ m = 1$ and $\beta_1 = 1$ | exponentiated generalized defined by |
| | | $\mu_1(.)(x) = (1-\theta)\prod_{j=1}^{m} G_j^{\beta_j}(x) + \theta$ | | Mudholkar *et al.* (1995) |
| **9S1C1.2** | $bt^{b-1}$ | $m_1(.)(x) = \theta\prod_{j=1}^{m} G_j^{\beta_j}(x)$ <br> $v_1(.)(x) = (1-\theta)\prod_{i=1}^{m}\left(1-G_i^{\alpha_i}(x)\right)^{\delta_i} + \theta$ | $\theta = 1,\ m = 1$ and $\beta_1 = 1$ | exponentiated generalized defined by Mudholkar *et al.* (1995) |
| **3S1C1.2** | $bt^{b-1}$ | $\ell_1(.)(x) = \theta\prod_{i=1}^{m}\left(1-G_i^{\alpha_i}(x)\right)^{\delta_i}$ <br> $\mu_1(.)(x) = (1-\theta)\prod_{j=1}^{m} G_j^{\beta_j}(x) + \theta$ | $\theta = 0$ and $m = 1$ | exponentiated generalized defined by Mudholkar *et al.* (1995) |
| **5S1C1.2** | ------------------ | $\mu_1(.)(x) = \prod_{i=1}^{m}\left(b_i + G_i^{\alpha_i}(x)\right)^{\beta_i}$ | $\theta = 0$ and $\beta = 1$ | exponentiated generalized defined by Mudholkar *et al.* (1995) |
| **3S1C1.2** | $abt^{a-1}(1-t^a)^{b-1}$ | $\ell_1(.)(x) = \theta\prod_{i=1}^{m}\left(1-G_i^{\alpha_i}(x)\right)^{\delta_i}$ <br> $\mu_1(.)(x) = (1-\theta)\prod_{j=1}^{m} G_j^{\beta_j}(x) + \theta$ | $\theta = 0,\ m = 1$ and $\beta_1 = 1$ | Kumaraswamy generalized defined by Cordeiro and Castro (2011) |
| **9S1C1.2** | $abt^{a-1}(1-t^{a})^{b-1}$ | $m_1(.)(x) = \theta \prod_{j=1}^{m} G_j^{\beta_j}(x)$<br>$v_1(.)(x) = (1-\theta) \prod_{i=1}^{m} \left(1 - G_i^{\alpha_i}(x)\right)^{\delta_i} + \theta$ | $\theta = 1,\ m = 1$ and $\beta_1 = 1$ | Kumaraswamy generalized defined by Cordeiro and Castro (2011) |
| **6S1C1.2** | ------------------ | $\ell_1(.)(x) = \prod_{i=1}^{m} \left(b_i - G_i^{\alpha_i}(x)\right)^{\beta_i}$ | $m = 1,\ b_1 = 1,\ \beta_1 = \beta$ and $\alpha_1 = \alpha$ | Kumaraswamy generalized defined by Cordeiro and Castro (2011) |
| **3S1C1.2** | $\frac{t^{a-1}(1-t)^{b-1}}{B(a,b)(1+t)^{a+b}}$ | $\ell_1(.)(x) = \theta \prod_{i=1}^{m} \left(1 - G_i^{\alpha_i}(x)\right)^{\delta_i}$<br>$\mu_1(.)(x) = (1-\theta) \prod_{j=1}^{m} G_j^{\beta_j}(x) + \theta$ | $\theta = 0,\ m = 1$ and $\beta_1 = 1$ | beta3 generalized defined by Thair and Nadarajah (2013) |
| | | | $\theta = 0$ and $m = 1$ | Mc3 generalized defined by Thair and Nadarajah (2013) |
| **9S1C1.2** | $\frac{t^{a-1}(1-t)^{b-1}}{B(a,b)(1+t)^{a+b}}$ | $m_1(.)(x) = \theta \prod_{j=1}^{m} G_j^{\beta_j}(x)$ <br> $\upsilon_1(.)(x) = (1-\theta)\prod_{i=1}^{m}\left(1-G_i^{\alpha_i}(x)\right)^{\delta_i} + \theta$ | $\theta = 1,\ m = 1$ and $\beta_1 = 1$ | beta1 generalized defined by Eugene *et al.* (2002) |
| | | | $\theta = 1$ and $m = 1$ | Mc3 generalized defined by Thair and Nadarajah (2013) |
| **3S1C1.2** | $\frac{t^{a-1}(1-t)^{b-1}exp(-ct)}{B(a,b)}$ | $\ell_1(.)(x) = \boldsymbol{\theta}\prod_{i=1}^{m}\left(1-G_i^{\alpha_i}(x)\right)^{\delta_i}$ <br> $\mu_1(.)(x) = (\mathbf{1}-\boldsymbol{\theta})\prod_{j=1}^{m} G_j^{\beta_j}(x) + \boldsymbol{\theta}$ | $\theta = 0,\ m = 1$ and $\beta_1 = 1$ | Kummer beta generalized defined by Pescim *et al.* (2012) |
| **9S1C1.2** | $\frac{t^{a-1}(1-t)^{b-1}exp(-ct)}{B(a,b)}$ | $m_1(.)(x) = \theta \prod_{j=1}^{m} G_j^{\beta_j}(x)$ <br> $\upsilon_1(.)(x) = (1-\theta)\prod_{i=1}^{m}\left(1-G_i^{\alpha_i}(x)\right)^{\delta_i} + \theta$ | $\theta = 1,\ m = 1$ and $\beta_1 = 1$ | Kummer beta generalized defined by Pescim *et al.* (2012) |
| **3S1C1.2** | $\frac{b^a}{\Gamma(a)} t^{a-1}e^{-bt}$ | $\ell_1(.)(x) = \theta\prod_{j=1}^{m}\left(1-G_j^{\beta_j}(x)\right)^{\gamma_j}$ | $\theta = 0,\ m = 1,\ \alpha_1 = 1,\ \lambda = 1,\ r = 1$ and $\delta_1 = 1$ | Generalized gamma-generated defined by Zografos (2009) |
| | | $\mu_1(.)(x) = \theta + \left(-ln\left(\prod_{i=1}^{m}\left(1 - G_i^{\alpha_i}(x)\right)^{\delta_i}\right)^{\lambda}\right)^{r}$ | | |
| **3S1C1.2** | $\frac{b^a}{\Gamma(a)} t^{a-1} e^{-bt}$ | $\ell_1(.)(x) = \theta \prod_{j=1}^{m}\left(1 - G_j^{\beta_j}(x)\right)^{\gamma_j}$ <br> $\mu_1(.)(x) = \theta + \left(-ln\left(1 - \prod_{i=1}^{m} G_i^{\alpha_i}(x)\right)^{\lambda}\right)^{r}$ | $\theta = 0$, $m = 1$, $\alpha_1 = 1$, $\lambda = 1$ and $r = 1$ | Generalized gamma-generated defined by Zografos (2009) |
| **3S1C1.2** | $\frac{b^a}{\Gamma(a)} t^{a-1} e^{-bt}$ | $\ell_1(.)(x) = \rho \prod_{i=1}^{m}\left(1 - G_i^{\omega_i}(x)\right)^{s_i}$ <br> $\mu_1(.)(x) = \rho - ln\left(1 - \prod_{l=1}^{m} G_l^{\lambda_l}(x)\right)$ | $\rho = 0, m = 1$, $\alpha = 0$ and $\lambda_1 = 1$ | Generalized gamma-generated defined by Zografos (2009) |
| **3S1C1.2** | $\frac{b^a}{\Gamma(a)} t^{a-1} e^{-bt}$ | $\ell_1(.)(x) = \rho \prod_{i=1}^{m}\left(1 - G_i^{\lambda_i}(x)\right)^{s_i}$ <br> $\mu_1(.)(x) = \rho + \left(-ln\left(\prod_{i=1}^{m}\left(1 - G_i^{\alpha_i}(x)\right)^{\omega_i}\right)^{\lambda}\right)^{r}$ | $m = 1$, $\alpha_1 = 1$, $\rho = 0$, $\omega_1 = 1$, $\lambda =$ 1 and $r = 1$ | Generalized gamma-generated defined by Zografos (2009) |
| **9S1C1.2** | $\frac{b^a}{\Gamma(a)} t^{a-1} e^{-bt}$ | $m_1(.)(x) = \rho \prod_{i=1}^{m} \left(1 - G_i^{\omega_i}(x)\right)^{s_i}$ <br> $v_1(.)(x) = \rho - ln\left(1 - \prod_{l=1}^{m} G_l^{\lambda_l}(x)\right)$ | $\rho = 0$, $m = 1$ and $\lambda_1 = 1$ | Generalized gamma-generated defined by Zografos (2009) |
| **9S1C1.2** | $\frac{b^a}{\Gamma(a)} t^{a-1} e^{-bt}$ | $m_1(.)(x) = \rho \prod_{i=1}^{m} \left(1 - G_i^{\lambda_i}(x)\right)^{s_i}$ <br> $v_1(.)(x) = \rho + \left(-ln\left(\prod_{i=1}^{m}\left(1 - G_i^{\alpha_i}(x)\right)^{\beta_i}\right)^{\lambda}\right)^{r}$ | $m = 1$, $\alpha_1 = 0$, $\rho = 0$ and $\beta_1 = 1$ | Generalized gamma-generated defined by Zografos (2009) |
| **9S1C1.2** | $\frac{b^a}{\Gamma(a)} t^{a-1} e^{-bt}$ | $m_1(.)(x) = \theta \prod_{j=1}^{m} \left(1 - G_j^{\beta_j}(x)\right)^{\gamma_j}$ <br> $v_1(.)(x) = \theta + \left(-ln\left(\prod_{i=1}^{m} G_i^{\alpha_i}(x)\right)\right)^{\delta}$ | $\theta = 0$, $m = 1$, $\alpha_1 = 1$ and $\delta = 1$ | Gamma G defined by Silva (2013) |
| **9S1C1.2** | $\frac{b^a}{\Gamma(a)} t^{a-1} e^{-bt}$ | $m_1(.)(x) = \theta \left(1 - \prod_{j=1}^{m} G_j^{\alpha_j}(x)\right)^{\lambda}$ <br> $v_1(.)(x) = \theta + \left(-ln\left(1 - \prod_{i=1}^{m}\left(1 - G_i^{\beta_i}(x)\right)^{\gamma_i}\right)^{r}\right)^{s}$ | $\theta = 0, m = 1$, $\beta_1 = 1, \gamma_1 = 1, r = 1$ and $s = 1$ | Gamma G defined by Silva (2013) |
| **3S1C1.2** | $\frac{b^a}{\Gamma(a)} t^{a-1} e^{-bt}$ | $\ell_1(.)(x) = \rho \prod_{i=1}^{m} G_i^{\omega_i}(x)$<br>$\mu_1(.)(x) = \rho - ln\left(\prod_{l=1}^{m} G_l^{\lambda_l}(x)\right)$ | $\rho = 0, m = 1$ and $\lambda_1 = 1$ | Gamma G defined by Silva (2013) |
| **3S1C1.2** | $\frac{b^a}{\Gamma(a)} t^{a-1} e^{-bt}$ | $\ell_1(.)(x) = \rho \prod_{i=1}^{m} G_i^{\lambda_i}(x)$<br>$\mu_1(.)(x) = \rho + \left(-ln\left(1 - \prod_{i=1}^{m}\left(1 - G_i^{\alpha_i}(x)\right)^{\omega_i}\right)^{\lambda}\right)^{r}$ | $m = 1$, $\rho = 0$, $\alpha_1 = 1$, $\delta_1 = 1$, $\lambda = 1$ and $r = 1$ | Gamma G defined by Silva (2013) |
| **9S1C1.2** | $\frac{b^a}{\Gamma(a)} t^{a-1} e^{-bt}$ | $m_1(.)(x) = \rho \prod_{i=1}^{m} G_i^{\omega_i}(x)$<br>$v_1(.)(x) = \rho - ln\left(\prod_{l=1}^{m} G_l^{\lambda_l}(x)\right)$ | $\rho = 0$, $m = 1$, $\alpha = 0$ and $\lambda_1 = 1$ | Gamma G defined by Silva (2013) |
| **9S1C1.2** | $\frac{b^a}{\Gamma(a)} t^{a-1} e^{-bt}$ | $m_1(.)(x) = \rho \prod_{i=1}^{m} G_i^{\lambda_i}(x)$<br>$v_1(.)(x) = \rho + \left(-ln\left(1 - \prod_{i=1}^{m}\left(1 - G_i^{\alpha_i}(x)\right)^{\omega_i}\right)^{\lambda}\right)^{r}$ | $m = 1$, $\rho = 0$, $\alpha_1 = 1$, $\omega_1 = 1$, $\lambda = 1$ and $r = 1$ | Gamma G defined by Silva (2013) |
| **12S1C1.2** | ------------------ | $m_1(.)(x) = \left( \frac{b\left(1 - G_2^{\beta}(x)\right)^{\gamma}}{G_1^{\alpha}(x) + b\left(1 - G_2^{\beta}(x)\right)^{\gamma}} \right)^{\theta}$ | $G_1(x) = G_2(x)$, $\alpha = 1$, $\beta = 1$ and $\theta = 1$ | Marshall-Olkin defined by Marshall and Olkin (1997) |
| | | | $G_1(x) = G_2(x)$, $\alpha = 1$ and $\beta = 1$ | Marshall-Olkin $G_1$ defined by Jayakumar and Mathew (2008) |
| **5S1C1.2** | ------------------ | $\mu_1(.)(x) = \left( \frac{G_1^{\alpha}(x)}{G_1^{\alpha}(x) + b\left(1 - G_2^{\beta}(x)\right)^{\gamma}} \right)^{\theta}$ | $G_1(x) = G_2(x)$, $\alpha = 1$, $\beta = 1$ and $\theta = 1$ | Marshall-Olkin defined by Marshall and Olkin (1997) |
| | | | $G_1(x) = G_2(x)$, $\alpha = 1$ and $\beta = 1$ | |
| | | | | Marshall-Olkin $G_1$ defined by Thair and Nadarajah (2013) |
| **12S1C1.2** | ------------------ | $m_1(.)(x) = exp\left(-\lambda \prod_{l=1}^{m} G_l^{\alpha_l}(x)\right)$ | $m = 1$ and $\alpha_1 = 1$ | Kumaraswamy $G_1$ Poisson defined by Ramos (2014) |
| **12S1C1.2** | ------------------ | $m_1(.)(x) = exp\left(-\lambda - \lambda \prod_{l=1}^{m} \left(1 - G_l^{\alpha_l}(x)\right)^{\beta_l}\right)$ | $m = 1$, $\alpha_1 = 1$ and $\beta_1 = 1$ | Kumaraswamy $G_1$ Poisson defined by Ramos (2014) |
| **12S1C1.2** | ------------------ | $m_1(.)(x) = exp\left(-\lambda - \lambda \left(1 - \prod_{l=1}^{m} G_l^{\alpha_l}(x)\right)^{\beta}\right)$ | $m = 1$, $\alpha_1 = 1$ and $\beta = 1$ | Kumaraswamy $G_1$ Poisson defined by Ramos (2014) |
| **2S1C1.2** | $\gamma t^{\gamma-1}$ | $u_1(.)(x) = \frac{e^{\lambda e^{-\beta x^{\alpha}}} - e^{\lambda}}{1-e^{\lambda}}$<br>$u_2(.)(x) = \frac{e^{-\frac{\lambda}{\alpha}W(-\alpha e^{-\alpha})} - e^{-\frac{\lambda}{\alpha}W(\psi(x))}}{e^{-\frac{\lambda}{\alpha}W(-\alpha e^{-\alpha})}-1}$,<br>$W(x) = \sum_{n=1}^{\infty} \frac{(-1)^{n-1} n^{n-2}}{(n-1)!} x^n$<br>$\psi(x) = -\alpha e^{-\alpha - bx^{a}}$<br>$\ell_1(.)(x) = \theta\left(1 - \frac{(1-\beta)^{-s} - \{1-\beta[1-G(x)]\}^{-s}}{(1-\beta)^{-s}-1}\right)$<br>$\mu_1(.)(x) = (1-\theta)\left(\frac{\zeta(s) - Li_s[1-G(x)]}{\zeta(s)}\right)^{\delta} + \theta,$<br>$Li_s(z) = \sum_{j=1}^{\infty} \frac{z^j}{j^s}$<br>$\zeta(s) = \sum_{j=1}^{\infty} \frac{1}{j^s}$ | $k=2$, $\alpha_1 = 1$, $\alpha_2 = 0$, $\theta = 0$ and $\delta = 0$ | Beta Weibull Poisson Family defined by Paixão (2014) |
| | | | $k=2$, $\alpha_1 = 0$, $\alpha_2 = 1$, $\theta = 0$ and $\delta = 0$ | Weibull Generalized Poisson class defined by Paixão (2014) |
| | | | $k=2$, $\alpha_1 = 0$, $\alpha_2 = 0$ and $\theta = 1$ | G-Negative Binomial family defined by Paixão (2014) |
| | | | $k=2$, $\alpha_1 = 0$, $\alpha_2 = 0$ and $\theta = 0$ | Zeta-G class defined by Paixão (2014) |
| **2S1C1.2** | $bt^{b-1}$ | $u_1(.)(x) = \sum_{j=0}^{x} \frac{C^{(j)}(a)}{j!\,C(\lambda)} (\lambda - a)^j$<br>$u_2(.)(x) = \sum_{j=1}^{x} \frac{1}{j!}\left[\left(C(0)\right)^j\right]^{(j-1)}$ | $k=2$, $\alpha_1 = 1$, $\alpha_2 = 0$, $\theta = 0$ and $\delta = 0$ | Power Series defined by Consul and Famoye (2006) |
| | | | $k=2$, $\alpha_1 = 0$, $\alpha_2 = 1$, $\theta = 0$ and $\delta = 0$ | Basic Lagrangian defined by Consul and Famoye (2006)) |
| | | $\ell_1(.)(x) = \theta\left(1 - \sum_{j=n}^{x} \frac{n}{(j-n)!\,j}\left[\left(C(0)\right)^j\right]^{(j-n)}\right)$ $\mu_1(.)(x) = (1-\theta)\left(\sum_{j=0}^{x} P(X=j)\right)^{\delta} + \theta$ $P(X=j) = \begin{cases} w(0), & j=0 \\ \left[\left(C(0)\right)^j w^{(1)}(0)\right]^{(j-1)}, & j=1,2,3,\ldots \end{cases}$ | $k=2$, $\alpha_1=0$, $\alpha_2=0$ and $\theta=1$ | Lagrangian Delta defined by Consul and Famoye (2006) |
| | | | $k=2$, $\alpha_1=0$, $\alpha_2=0$ and $\theta=0$ | Generalized Lagrangian family (CONSUL and FAMOYE, 2006) |
| **2S1C1.2** | $bt^{b-1}$ | $u_1(.)(x) = \int_{-\infty}^{x} e^{\int \frac{a_0+a_1t+\cdots+a_st^s}{b_0+b_1t+\cdots+b_rt^r} f(t)dt}\, dt$ $\ell_1(.)(x) = \theta\left(1 - \int_{-\infty}^{x} e^{\int \frac{a_0+a_1t+\cdots+a_st^s}{b_0+b_1t+\cdots+b_rt^r} dt}\, dt\right)$ $\mu_1(.)(x) = (1-\theta)\left(\int_{-\infty}^{x}\int_{-\infty}^{y}\left(\sum_{i=1}^{2}\alpha_i(t) f^{\beta_i}(t)\right)dt\,dy\right)^{\delta} + \theta$ | $k=1$, $\alpha_1=1$, $\theta=0$ and $\delta=0$ | Generalized Pearson in Ordinary Differential Equation form defined by Shakil *et al.* (2010) |
| | | | $k=1$, $\alpha_1=0$ and $\theta=1$ | Generalized Pearson in Ordinary Differential Equation form defined by Shakil *et al.* (2010) |
| | | | $k=1$, $\alpha_1=0$, $\theta=0$ and $\delta=1$ | Generalized Family in Ordinary Differential |
| | | | | Equation form defined by (VODA, 2009) |